\documentclass[12pt]{amsart}
\usepackage{amssymb,latexsym, amsthm,epsf}

\newtheorem{thm}{Theorem}[section] 
\newtheorem{defn}[thm]{Definition}

\newtheorem{cor}[thm]{Corollary}

\newtheorem{lem}[thm]{Lemma}
\newtheorem{exa}[thm]{Example}
\newtheorem{prop}[thm]{Proposition}
\newtheorem{rem}[thm]{Remark}
\newtheorem{algo}[thm]{Algorithm}

\def\R{{\mathbb R}}
\def\P{{\mathbb P}}

\def\Z{{\mathbb Z}}

\def\({\left(}
\def\){\right)}

\long\def\forget#1\forgotten{}

\newif \iffurther 
\furtherfalse

\newif \iffurther 
\furtherfalse

\newif\ifXY 
\XYfalse    

\ifXY
\usepackage{xy}
\fi
\ifXY
\xyoption{all}
\fi

\begin{document}

\title{Spindle-configurations of skew lines}

\author[Roland Bacher and David Garber]{Roland
  Bacher$^{1,2}$ and David Garber$^3$}

\stepcounter{footnote} 
\footnotetext{Corresponding author.} 
\stepcounter{footnote} 
\footnotetext{Support from the Swiss National Science Foundation is 
gratefully acknowledged.}
\stepcounter{footnote} 
\footnotetext{Partially supported by the Chateaubriand postdoctoral
  fellowship funded by the French government.} 

\address{Institut Fourier, BP 74, 38402 Saint-Martin D'Heres CEDEX,
  France}
\email{\{bacher,garber\}@mozart.ujf-grenoble.fr}

\keywords{Configurations of skew lines, Spindles, linking matrix,
  switching graphs, permutation}

\begin{abstract} 
We prove a conjecture of Crapo and Penne which characterizes isotopy
classes of skew configurations with spindle-structure. We use this
result in order to define an invariant, spindle-genus, for 
spindle-configurations.

We also slightly simplify the exposition of some known invariants for 
configurations of skew lines and use them to define a natural
partition of the lines in a skew configuration. 

Finally, we describe an algorithm which constructs a spindle 
in a given switching class, or proves non-existence of such a spindle. 
\end{abstract}

\maketitle

\section{Introduction}

{\it A configuration of $n$ skew lines in $\R^3$} or a {\it skew 
configuration} is a set
of $n$ non-intersecting lines in $\R^3$
containing no pair of parallel lines.

Two skew configurations $C_1$ and $C_2$ are 
{\it isotopic} if there exists an isotopy (continuous deformation of 
skew configurations) from $C_1$ to $C_2$. 

The study and classification of configurations of skew lines 
(up to isotopy) was started by Viro \cite{V} and continued for example in
\cite{BM1}, \cite{BM2}, \cite{CP}, 
\cite{K}, \cite{MSu}, 
\cite{P1}, \cite{P2}, \cite{P3}, \cite{P4} and \cite{VD}.
 
A {\it spindle} (or {\it isotopy join} or {\it horizontal configuration})
is a particularly nice configuration of skew lines 
in which all lines intersect an oriented additional 
line $A$, called the {\it axis} of the spindle. 
Its isotopy class is completely described by a
{\it spindle-permutation} $\sigma:\{1,\dots,n\}\longrightarrow \{1,\dots,n\}$
encoding the order in which an open half-plane 
revolving around its boundary $A$ intersects the lines 
during a half-turn (see Section 
\ref{spindles} for the precise definition).
A {\it spindle-configuration} is a skew configuration isotopic to a spindle.

Consider the {\it spindle-equivalence} relation on permutations 
of $\{1,\dots,n\}$
generated by transformations of the following three types:

\begin{enumerate}
\item $\sigma\sim \mu$ if $\mu(i)=s+\sigma(i+t\pmod n) \pmod n$ for
some integers $0\leq s,t<n$ (all sums are modulo $n$).

\item $\sigma\sim \mu$ if $\sigma(i)\leq k$ for $i\leq k$ and
$$\mu(i)=\left\{\begin{array}{ll}
k+1-\sigma(k+1-i)\qquad &i\leq k\cr
\sigma(i)&i> k
\end{array}\right.$$
for some integer $k\leq n$.

\item  $\sigma\sim \mu$ if $\sigma(i)\leq k$ for $i\leq k$ and
$$\mu(i)=\left\{\begin{array}{ll}
\sigma^{-1}(i)\qquad &i\leq k\cr
\sigma(i)&i> k\end{array}\right.$$
for some integer $k\leq n$.
\end{enumerate}

Conjecture 59 in \cite{CP} states that two spindle-configurations 
are isotopic if and only if they 
are described by spindle-equivalent permutations. Its proof
is the main result of this paper:

\begin{thm} \label{main}
Two spindle-permutations $\sigma,\sigma'$ 
give rise to isotopic spindle-configurations
if and only if $\sigma$ and $\sigma'$ are spindle-equivalent.
\end{thm}

Orienting and labeling all lines of a skew configuration, one gets a {\it
linking matrix} by considering the signs of crossing lines.
The associated {\it switching class} or {\it homology equivalence class}
is independent of labelings and orientations. 
A result of Khashin and Mazurovskii, Theorem 3.2 in \cite{KM},
states that homology-equivalent spindles (spindles defining
the same switching class) are isotopic. We have thus:

\begin{cor} \label{maincor}
Two spindle permutations define the same switching class if
and only if they are spindle-equivalent.
\end{cor}

Isotopy classes of spindle-configurations have thus an easy 
combinatorial description and can be considered as \lq\lq understood'',
either in terms of spindle-equivalence classes or in terms of
switching classes,
in contrast to the general case where no (provenly) complete invariants 
are available.
\medskip

An {\it invariant} is a map
$$\{\hbox{Configurations of skew lines up to isotopy}\}
\longrightarrow {\mathcal R}$$
where ${\mathcal R}$ is an algebraic (or combinatorial) structure, 
usually a group, vector-space or ring (or a finite set). An invariant is
{\it complete} if the corresponding map is injective.

A few useful invariants for configurations of skew lines are:
\begin{enumerate}
 
\item {\it Equivalence classes of skew pseudoline diagrams} 
(see Section \ref{diagrams}): Completeness
unknown (this is a major problem
since the obvious planar representation of skew configurations
is perhaps not faithful).   
A powerful combinatorial invariant somewhat tedious to 
handle. Switching classes and Kauffman polynomials factorize
through it.

\item {\it Switching classes} or {\it two-graphs}, see page 7 of \cite{Za1},
also called {\it homological equivalence classes}, see
\cite{BM2}, are equivalent to
the description of the sets of {\it linking numbers} or 
{\it linking coefficients}, see \cite{BM2} or \cite{VD}.
The definition of this invariant uses a
{\it linking matrix} encoding the signs of oriented crossing 
(not intersecting in the compactification
$\R\P^3\supset \R^3$) lines. 
The switching class is a complete invariant for configurations
of up to $5$ skew lines and is not complete for more than $5$ lines.
Theorem 3.2 of \cite{KM} states that switching-classes are a complete 
invariant for spindle-configurations.

Many weaker invariants can be derived from switching classes: Section 
\ref{euler} describes 
the {\it Euler partition} of a switching class which yields a natural
set-partition on the lines of a skew configuration. The definition
of the Euler partition is more natural (and better known) for an odd
number of lines.

Slightly weaker (but more elementary to handle) than the switching
class is the characteristic polynomial of a linking matrix $X$
$$P_X(t)=\sum_{i=0}^n\alpha_i t^i=\hbox{det}(t {\rm I}-X)\ .$$
The spectrum $\hbox{spec}(X)=\{\lambda_1\leq \lambda_2\leq\dots\leq 
\lambda_n\}$ of $X$ or the traces $\hbox{tr}(X)=\sum_{i=1}^n \lambda_i^k,\ 
k=1,\dots,n$ of its first powers yield of course the same information.

The coefficient $\alpha_{n-3}$ of $P_X(t)$
conveys the same information as {\it chirality},
a fairly weak invariant considered for instance in
\cite[Section 3]{CP}.

\item {\it Kauffman polynomials}: 
Completeness unknown. A powerful invariant which is
unfortunately difficult to compute, see \cite{D} or 
\cite[Section 14 and Appendix]{CP}
for the definition and examples.

\item {\it Link invariants} for links in the $3-$sphere
${\mathbb S}^3$ applied to the
preimage $\pi^{-1}(C)\subset {\mathbb S}^3$ (called a {\it Temari model} by
some authors, see for instance Section 11 of \cite{CP}) 
of a skew configuration $C\subset \R^3\subset \R\P^3$ under
the double covering $\pi:{\mathbb S}^3\longrightarrow \R\P^3$.

\item Existence of a {\it spindle structure}. A generally
very weak invariant since spindle structures are rare among configurations
with many lines. Theorem \ref{main} shows however that the spindle-equivalence
class provides a complete invariant for the very small
subset of skew configurations with a spindle structure.

\item One should also mention the {\it shellability order}, a 
generalization of the notion 
of spindle structure, used as a classification-tool in \cite{BM2}.

\end{enumerate}
 
\medskip

The sequel of the paper is organized as follows:

Section \ref{diagrams} introduces skew pseudoline diagrams.

Sections \ref{mat_invar}-\ref{euler} are devoted to (various aspects of) 
switching classes.

Section \ref{spindles} describes spindle-configurations and 
contains a proof of the easy (and known) direction in
Theorem \ref{main}: Spindle-equivalent permutations yield 
isotopic configurations.

Section \ref{isospindles} proves the difficult direction: Isotopic
spindles have spindle-equivalent permutations. This completes the proof
of Theorem \ref{main}.

Section \ref{KaMa} recalls (for self-containedness) 
the proof of Theorem 3.2 in \cite{KM}:
Spindle configurations with switching-equivalent linking matrices are
isotopic. Corollary \ref{maincor} follows.

Section \ref{spindlegenus} describes a somewhat curious invariant for spindle-equivalence
classes of permutations (or spindle-configurations) which involves 
$2-$dimensional topology.

In Section \ref{algo} we describe a fast algorithm which computes (or proves
non-existence of) a spindle-permutation (which is unique up to 
spindle-equivalence by Corollary \ref{maincor}) having a linking matrix of
given switching class.   

Section \ref{enumeration} contains a few computational data.


\section{Skew pseudoline diagrams}\label{diagrams}
Skew pseudoline diagrams are combinatorial objects
providing a convenient tool for studying configurations of skew lines.

\begin{defn}
A {\bf pseudoline} in $\R^2$ is a smooth simple
curve representing a non-trivial cycle in $\R \P^2$.
An {\bf arrangement of $n$ pseudolines} in $\R^2$ is a set of $n$
pseudolines with pairs of pseudolines intersecting
transversally exactly once. An arrangement is {\bf generic} if no
triple intersections occur.
\end{defn}

\begin{defn}
{\bf A skew pseudoline diagram of $n$ pseudolines} in $\R^2$ is a generic
arrangement of $n$ pseudolines with crossing data at intersections. 
The crossing data selects at each intersection the over-crossing
pseudoline.
\end{defn}
 
We draw skew pseudoline diagrams with the conventions 
used for knots and links: under-crossing curves are
slightly interrupted at crossings.

Skew pseudoline diagrams are {\it equivalent} if they are related
by a finite sequence of the following moves (see \cite[Section 9]{CP}):
\begin{enumerate}
\item Reidemeister-3  (or $*$-move), the most interesting of the
three classical moves for knots and links (see Figure \ref{R-S3}).

\begin{figure}[h]\label{R-S3}
\epsfysize=2cm
\epsfbox{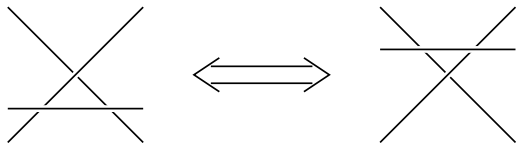}
\caption{Local description of a Reidemeister-3 move}
\end{figure}                                                                 

\item Projective move (or $||$-move): pushing a crossing through
  infinity
(see Figure \ref{infinity}).
\begin{figure}[h]\label{infinity}
\epsfysize=2.5cm
\epsfbox{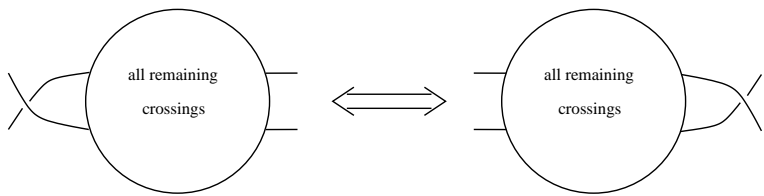}
\caption{Projective move}
\end{figure}                                                                 

\end{enumerate}

Generic projections of isotopic skew configurations yield
equivalent skew pseudoline diagrams (cf. for instance \cite[Theorem 48]{CP}).  

Not every (equivalence class of a) skew pseudoline diagram arises
by projecting of a suitable skew configuration: a configuration
involving at least $4$ skew lines never projects on a 
skew pseudoline diagram which is alternating (see \cite{PPW}). 
There are even generic arrangements of $n\geq 9$
pseudolines which are not stretchable, i.e. which cannot be realized 
as an arrangement of straight lines, see \cite{Gr} for an example with 9 
pseudolines.

The existence of non-isotopic skew configurations
inducing equivalent skew pseudoline diagrams is unknown, cf.
\cite[Section 17, Problem 2]{CP}. 
This is a major difficulty for classifications: All known invariants for
skew configurations factor through skew pseudoline diagrams.

\section{Linking matrices and switching classes}\label{mat_invar}

\begin{figure}[h]\label{fig1}
\epsfysize=3cm
\epsfbox{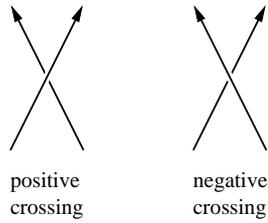}
\caption{Positive and negative crossings}
\end{figure}

One assigns signs to pairs of oriented under- or 
over-crossing curves (as arising for instance from oriented
knots and links) which are drawn in the oriented plane. 
Figure \ref{fig1} shows a positive and a negative crossing.

The {\it sign} or {\it linking number} 
$\hbox{lk}(L_A,L_B)$ between two oriented skew lines 
$L_A,L_B\subset \R^3$
can be computed as follows: choose ordered pairs of points 
$(A_\alpha,A_\omega)$ on $L_A$ (resp.
$(B_\alpha,B_\omega)$ on $L_B$) which induce the orientations. The sign of the
crossing determined by $L_A$ and $L_B$ is then given by
$$\hbox{lk}(L_A,L_B)=\hbox{sign}\Big(\hbox{det}\left(\begin{array}{c}
A_\omega-A_\alpha\cr
B_\alpha-A_\omega\cr
B_\omega-B_\alpha\end{array}\right)\Big)\in \{\pm 1\}$$
where $\hbox{sign}(x)=\frac{x}{\vert x\vert}$ for  $x\not=0$.
 
Signs of crossings are also defined in skew pseudoline diagrams. 

The {\it linking matrix} of a diagram of $n$
oriented and labeled skew pseudolines
$L_1,\dots,L_n$ is the symmetric $n\times n$
matrix $X$ with diagonal coefficients $x_{i,i}=0$ and
$x_{i,j}=\hbox{lk}(L_i,L_j)$ for $i\not= j$.

\begin{figure}[h]
\epsfysize=5cm
\epsfbox{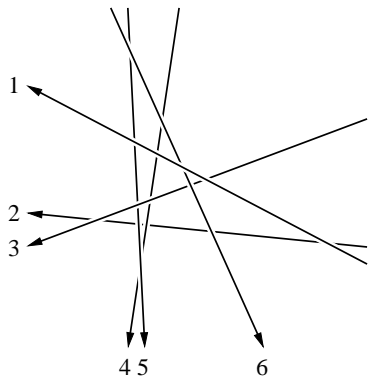}
\caption{A configuration of $6$ labeled and oriented skew lines}
\label{fig2}\end{figure}                                                               

Figure \ref{fig2} shows a labeled and oriented configuration of six skew
lines with linking matrix
$$X=\left( 
\begin{array}{rrrrrr}
0 &  1 &  1 &  1 &  1 &  1 \\  
1 &  0 & -1 & -1 & -1 & -1 \\  
1 & -1 &  0 &  1 &  1 & -1 \\  
1 & -1 &  1 &  0 & -1 & -1 \\  
1 & -1 &  1 & -1 &  0 & -1 \\  
1 & -1 & -1 & -1 & -1 &  0 
\end{array}
\right)\ .
$$

Two symmetric matrices $X$ and $Y$ are {\it switching-equivalent} if 
$$Y=D\ P^t\ X\ P\ D$$
where $P$ is a permutation matrix and $D$ is a diagonal matrix with
diagonal coefficients in $\{\pm 1\}$. Since $PD\in O(n)$ is
orthogonal, we have
$(PD)^{-1}=DP^t$. Switching-equivalent matrices are thus conjugate and
share a common characteristic polynomial.

\begin{prop} All linking matrices of a fixed skew pseudoline diagram
are switching-equivalent.
\end{prop}

\begin{proof}
Relabeling the lines conjugates a linking matrix $X$
by a permutation matrix. Reversing the orientation
of some lines amounts to conjugation by a diagonal $\pm 1$ matrix.
\end{proof}

\begin{rem}\label{switchingterminology} The terminology \lq\lq
switching classes'' (many authors use also \lq\lq two-graphs'')
is motivated by the following combinatorial interpretation and
definition of switching classes: 

Two finite simple (loopless and no multiple edges)
graphs $\Gamma_1$
and $\Gamma_2$ with vertices $V$ and (unoriented)
edges $E_1,E_2\subset V\times V$
are {\it switching-related} with respect to a subset $V_-\subset V$ of 
vertices if their edge-sets
$E_1,E_2$ coincide on $(V_-\times V_-)\cup \big((V\setminus V_-)\times
(V\setminus V_-)\big)$ and are complementary on $\big(V_-\times (V\setminus V_-)\big)\cup \big((V\setminus V_-)
\times V_-\big)$. A {\it switching class} of graphs 
is an equivalence class of switching-related graphs. 

\begin{figure}[h]
\epsfysize=2.5cm
\epsfbox{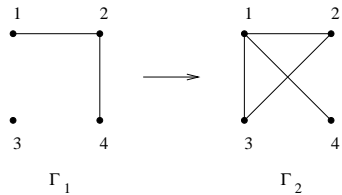}
\caption{$\Gamma_1$ and $\Gamma_2$ are switching-related with respect to 
$\{1,2\}\subset \{1,2,3,4\}$}\label{switching}
\end{figure}                                                               

Encoding adjacency, resp. non-adjacency, of distinct vertices
by $\pm 1$ yields a bijection between switching classes of graphs
and switching classes of matrices.
Conjugation by permutation-matrices corresponds to relabeling the vertices of 
a graph $\Gamma$ and conjugation by a diagonal $\pm 1-$matrix 
corresponds
to the substitution of $\Gamma$ by a switching-related graph.
\end{rem}

\begin{rem}\label{character} The characteristic polynomial of a
linking matrix of a skew pseudoline diagram 
(or of a configuration of skew lines) is in general weaker than its 
switching class: The linking matrices
$$\left(
\begin{array}{rrrrrrrr}
 0 &   1 &   1 &   1 &   1 &   1 &   1 &   1\\   
 1 &   0 &  -1 &   1 &  -1 &   1 &  -1 &   1\\   
 1 &  -1 &   0 &   1 &   1 &   1 &   1 &  -1\\   
 1 &   1 &   1 &   0 &  -1 &   1 &   1 &  -1\\   
 1 &  -1 &   1 &  -1 &   0 &   1 &   1 &  -1\\   
 1 &   1 &   1 &   1 &   1 &   0 &  -1 &  -1\\   
 1 &  -1 &   1 &   1 &   1 &  -1 &   0 &  -1\\   
 1 &   1 &  -1 &  -1 &  -1 &  -1 &  -1 &   0   
\end{array}\right)
$$
and 
$$\left(
\begin{array}{rrrrrrrr}
 0 &   1 &   1 &   1 &   1 &   1 &   1 &   1\\   
 1 &   0 &   1 &   1 &  -1 &   1 &  -1 &   1\\   
 1 &   1 &   0 &   1 &  -1 &   1 &  -1 &   1\\   
 1 &   1 &   1 &   0 &  -1 &   1 &   1 &  -1\\   
 1 &  -1 &  -1 &  -1 &   0 &   1 &   1 &  -1\\   
 1 &   1 &   1 &   1 &   1 &   0 &  -1 &  -1\\   
 1 &  -1 &  -1 &   1 &   1 &  -1 &   0 &  -1\\   
 1 &   1 &   1 &  -1 &  -1 &  -1 &  -1 &   0   
\end{array}\right)
$$
are in different switching classes (see Example \ref{exam_euler_even}),
but have the same characteristic polynomial
$$(t-3)(t-1)^2(t+1)(t+3)^2(t^2-2t-11)\ .$$

This example is minimal: Distinct switching
classes of order less than $8$ have distinct characteristic polynomials.
\end{rem}

\begin{rem}\label{rem_dif_config_same_switching} Isotopy classes of
skew configurations with  $\leq 5$ lines are 
characterized by their switching class.

However, Figure \ref{2_configs} shows two non-isotopic 
skew configurations of $6$ lines (their Kauffman polynomials, 
see \cite[Section 14 and Appendix]{CP},
$$5A^{12}B^3 + 10A^{11}B^4 - 10A^9 B^6 + A^8 B^7 + 16 A^7 B^8 + 10A^6 B^9$$
$$- 6A^5 B^{10} - 5A^4 B^{11} + 6A^3 B^{12} + 6A^2 B^{13} - B^{15}$$
and
$$-A^{15} + 6A^{13} B^2 + 6A^{12} B^3 - 5A^{11} B^4 - 6A^{10} B^5$$
$$+10A^9 B^6 + 16A^8 B^7 + A^7 B^8 - 10A^6 B^9 + 10 A^4 B^{11} + 
5A^3 B^{12}$$
are different) having a common switching class represented by
$$\left( 
\begin{array}{rrrrrr}
0 &  1 &  1 &  1 &  1 &  1 \\  
1 &  0 &  1 & -1 & -1 &  1 \\  
1 &  1 &  0 &  1 & -1 & -1 \\  
1 & -1 &  1 &  0 &  1 & -1 \\  
1 & -1 & -1 &  1 &  0 &  1 \\  
1 &  1 & -1 & -1 &  1 &  0 
\end{array}
\right)
$$

\begin{figure}[h]
\epsfysize=5cm
\epsfbox{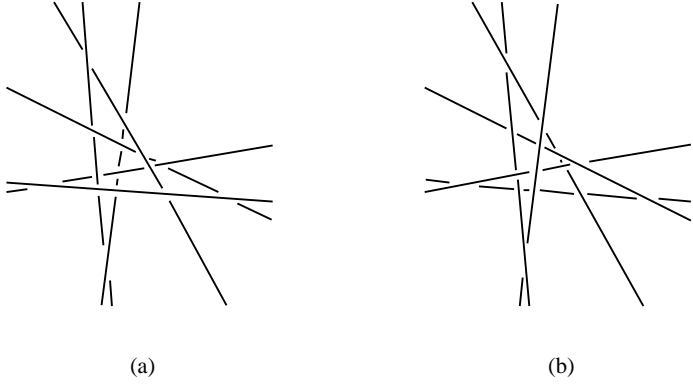}
\caption{Non-isotopic mirror configurations with switching-equivalent
linking matrices}\label{2_configs}
\end{figure}                                                               

\end{rem}

All $2^{n\choose 2}$ possible linking matrices can be obtained from a
fixed skew pseudoline diagram by an appropriate choice of crossing data. 
The number of equivalence classes of
skew pseudoline diagrams with $n$ pseudolines equals thus at least the 
number of switching classes of order $n$. 
Remark \ref{rem_dif_config_same_switching} shows that this 
inequality is in general strict.


\medskip

The mirror configuration $\overline C$ obtained by reflecting a 
skew configuration $C$ through the $z=0$ hyperplane 
has opposite crossing data. 
We have thus $\overline X=-X$ for the associated linking matrices.
An example of two such configurations $C,\overline C$ is also given by 
Figure \ref{2_configs}.

A configuration $C$ is {\it amphicheiral} if it is isotopic to its
mirror $\overline C$.

\begin{prop}\label{mirror} 
\begin{enumerate}
\item[(i)] The linking matrix $X$ of an amphicheiral
configuration of skew lines is switching-equivalent to $-X$. 
In particular, amphicheiral configurations containing an odd number 
of skew lines have non-invertible linking matrices.

\item[(ii)] Amphicheiral configurations with $n\equiv 3\pmod 4$ 
lines do not exist.
\end{enumerate}
\end{prop}

\begin{proof}
Assertion (i) is obvious.

Assertion (ii) is \cite[Theorem 1]{VD}. We rephrase the
proof using properties of linking matrices.

Let $\sum_{i=0}^n \alpha_i t^i=\det(tI-X)$ be the characteristic polynomial 
of the linking matrix $X$ for an amphicheiral configuration with $n$
skew lines. Assertion (i) shows that we have
$\alpha_{n-1}=\alpha_{n-3}=\alpha_{n-5}=\dots= 0$. This implies
$$0=\alpha_{n-3}=-\sum_{1\leq i,j,k\leq n}\left(x_{i,j}x_{j,k}x_{k,i}+x_{i,k}x_{k,j}x_{j,i}\right)=
-2\left(\sum_{1\leq i<j<k\leq n}x_{i,j}x_{j,k}x_{k,i}\right)\ .$$
For $n\equiv 3\pmod 4$\  the number $n\choose 3$ of summands in 
$\sum_{1\leq i<j<k\leq n}x_{i,j}x_{j,k}x_{k,i}$ is odd.
Since all these summands are $\pm 1$, we get a contradiction.
\end{proof}

\section{Switching classes and linking numbers}\label{linking}

{\it Linking numbers} (also called {\it homological equivalence 
classes or chiral signatures}) are a classical and
well-known invariant for skew pseudoline diagrams. We sketch below
briefly the well-known proof that they correspond to the switching 
class of a linking matrix.

In this paper we work with switching classes mainly 
because they are easier to handle. 

\medskip

The {\it linking number} $\hbox{lk}(L_i,L_j,L_k)$ of three lines 
(\cite{BM2},\cite{VD}) is defined as
the product $x_{i,j}x_{j,k}x_{k,i}\in \{\pm 1\}$ of the signs 
for the corresponding three crossings. The result 
is independent of the chosen orientations for $L_i,L_j$ and $L_k$ and
yields an invariant 
$$\{\hbox{triplets of lines in skew pseudoline diagrams}\}
\longrightarrow \{\pm 1\}\ .$$
  
{\it The set of linking numbers} is the list of the numbers 
$\hbox{lk}(L_i,L_j,L_k)$ for  all
triplets $\{L_i,L_j,L_k\}$ of lines in a skew pseudoline diagram. 

Linking numbers (defining a {\it two-graph}, see \cite{Za1}) and
switching classes are equivalent. Indeed, linking numbers of
a diagram $D$ of skew lines can easily be retrieved from 
a linking matrix for
$D$. Conversely, given all linking numbers ${\rm lk}(L_i,L_j,L_k)$
of a diagram $D$, choose an orientation of the first line $L_1$.
Orient the remaining lines 
$L_2,\dots,L_n$ such that they cross $L_1$ positively. 
A linking matrix $X$ for $D$ is then given by $x_{1,i}=x_{i,1}=1,\ 2\leq i\leq n$
and $x_{a,b}={\rm lk}(L_1,L_a,L_b)$ for $2\leq a\not= b\leq n$.

Two skew pseudoline diagrams are {\it homologically equivalent} if there
exists a bijection
between their lines, which preserves all linking numbers. 
Two diagrams are homologically equivalent if and only
if they have switching-equivalent linking matrices.

A last invariant considered by some authors (see \cite[Section 3 and
Appendix]{CP}) 
is the {\it chirality} $(\gamma_+,\gamma_-)$ of a skew pseudoline
diagram. It is defined 
as
$$\begin{array}{l}
\displaystyle \gamma_+=\sharp\{1\leq i<j<k\leq n\ \vert\ \hbox{lk}(L_i,L_j,L_k)=1\}
\cr
\displaystyle \gamma_-=\sharp\{1\leq i<j<k\leq n\ \vert\ \hbox{lk}(L_i,L_j,L_k)=-1\}
\end{array}$$
One has of course
$$\gamma_+=\frac{{n\choose 3}+c}{2}\ ,\quad  
\gamma_-=\frac{{n\choose 3}-c}{2}$$
where 
$$c=\sum_{1\leq i<j<k\leq n}x_{i,j}x_{j,k}x_{k,i}=
\frac{1}{6}\hbox{trace} (X^3)=-\frac{\alpha_{n-3}}{2}$$
is proportional to the coefficient of $t^{n-3}$ in the characteristic
polynomial $\det(tI-X)=\sum_{i=0}^n \alpha_it^i$ of a linking matrix $X$.

The sign indeterminancy in linking matrices representing switching classes
makes their use tedious. For switching classes of odd order,
a satisfactory answer addressing this problem will be given
in the next section. For even orders, there seems to be no
completely satisfactory way to get rid of all sign-indeterminancies.

One possible normalization consists of choosing a given (generally the
first) row of a representing matrix and to make all entries in this
row positive by conjugation with a suitable diagonal $\pm 1$ 
matrix. The resulting matrix describes a simple graph on
$n-1$ vertices encoding all entries outside the chosen row
(and column). The choice of the suppressed row leads to the notion
of graphs which are called \lq\lq cousins'' in \cite{CP} where this
point of view is adopted instead of switching-equivalence.

\section{Signature for configurations of $2n-1$ lines}

Consider the {\it signature}
$$\epsilon(X)=\prod_{1\leq i<j\leq 2n-1} x_{i,j}\in\{\pm 1\}$$
of a linking matrix $X$ of odd order $2n-1$. 
It is easy to check 
that $\epsilon(X)$ depends only on the switching-class of $X$. 
(This follows also from 
Proposition \ref{oddsignature} below.)

\begin{prop} \label{oddsignature} 
We have 
$$\epsilon(X)=(-1)^{\gamma_-}$$
where 
$\gamma_-=\sharp\{1\leq i<j<k\leq 2n+1\ \vert\ \hbox{lk}(L_i,L_j,L_k)=-1
\}$ counts the number of triplets of pseudolines with linking number
$-1$ in a pseudo-line diagram with linking matrix $X$.
\end{prop} 

\begin{proof} The identity is correct if $x_{i,j}=1$ for 
all $1\leq i<j\leq 2n-1$.
Reversing exactly one pair of coefficients $x_{i,j},x_{j,i}$ 
reverses the sign of the $2n-3$ linking numbers $lk(L_i,L_j,L_k,\ 
k\not= i,j$ and thus changes the parity of 
$\gamma_-$. This implies the result by induction on $\gamma_-$.
\end{proof}
 
\begin{rem} If the matrix $D$ corresponds to a configuration of $2n-1$ 
pseudolines arising from a spindle-configuration with spindle 
permutation $\sigma$, then $\epsilon(D)$ corresponds to the signature
of the permutation $\sigma$ (which is well-defined for the 
spindle-equivalence class of $\sigma$).
This fact provides another proof of assertion (ii) in 
Proposition \ref{mirror} for spindle-configurations: 
Given a spindle-permutation $\sigma$ of $\{1,\dots,2n-1\}$,
a spindle-permutation of the mirror configuration is for instance
given by $\sigma\tau$ where $\tau$ is the involution 
defined by $\tau(i)=2n-i,\ 1\leq i<2n$ and $\tau$ has 
signature $-1$ if $n\equiv 0\pmod 4$. 
\end{rem} 

\begin{rem} Defining invariants of switching classes is fairly easy.
A few examples are:
\begin{enumerate}

\item
The set of numbers $\vert\alpha_{i,j}\vert,\ i\not=j$ where 
$\alpha_{i,j}=\sum_{k=1}^n x_{i,k}x_{j,k}$ or the set of all triplets
$\alpha_{i,j}\alpha_{i,k}\alpha_{j,k}$.

\item
The set of all $4-$tuplets $\vert\sum_{k=1}^n x_{i_1,k}x_{i_2,k}
x_{i_3,k}x_{i_4,k}\vert$ with $i_1<i_2<i_3<i_4$.

\item
Properties of the (non-Euclidean) Lattice with Gram matrix 
(scalar products between generators) $X$ representing a given 
switching class.

\item
Properties of the Euclidean lattice spanned by the rows of $X$,
considered as integral vectors in standard Euclidean space.
\end{enumerate}
\end{rem}

\section{Euler partitions}\label{euler}

In this section we study invariants of switching classes which 
have a computational cost of $O(n^2)$ operations.
For skew configurations involving a huge number of lines they are thus
easier to compute and to compare
than chirality-invariants (with computational
cost $O(n^3)$).

The behaviour of switching classes depends on 
the parity of their order $n$.

Switching classes of odd order $2n-1$ are in bijection with Eulerian
graphs. This endows pseudoline diagrams consisting of an odd number
of pseudolines with a canonical orientation (up to a global change).
We get a partition of the pseudolines into equivalence
classes according to the number of positive crossings in which they
are involved for an Eulerian orientation.

We consider the case of odd order in Subsection \ref{oddswitching}.

The situation is more complicated for switching classes 
of even order $2n$. We replace Euler graphs appearing for odd orders 
by a suitable kind of planar rooted binary trees which we call 
{\it Euler trees}. The leaves of the Euler tree induce again a
natural partition, which we call the {\it Euler partition}, 
of the set of pseudolines into equivalence classes
of even cardinalities. Subsection 
\ref{even_case} addresses the even case.

\subsection{Switching classes of odd order - Euler orientations}
\label{oddswitching}

A simple finite graph $\Gamma$ is {\it Eulerian} if all its
vertices are of even degree. The following well-known result goes back
to Seidel \cite{Se}.

\begin{prop} \label{Eulerbijection}
Eulerian graphs with an odd number $2n-1$ of vertices are
in bijection with switching classes of order $2n-1$.
\end{prop}

We recall here the simple proof since it yields a fast 
algorithm for computing
Eulerian orientations on configurations with an odd number of skew lines.

\begin{proof}
Choose a representing matrix $X$ in a given switching class.
For $1\leq i\leq 2n-1$ define the number
$$v_i=\sharp\{j\ \vert\ x_{i,j}=1\}=\sum_{j=1,j\not= i}^{2n-1}
\frac{x_{i,j}+1}{2}$$
counting all entries equal to $1$ in the $i-$th row of $X$.
Since $X$ is symmetric, the vector $(v_1,\dots,v_{2n-1})$ 
has an even number of odd
coefficients and conjugation of the matrix $X$ with the diagonal matrix having 
diagonal entries $(-1)^{v_i}$ turns $X$ into a matrix $X_E$ with an even number
of $1$'s in each row and column. The matrix $X_E$ is well defined up 
to conjugation
by a permutation matrix and defines an Eulerian graph with vertices
$\{1,\dots,2n-1\}$ and edges $\{i,j\}$ if $(X_E)_{i,j}=1$. 
This construction can easily be reversed. 
\end{proof}

Let $D$ be a skew diagram having an odd number of pseudolines. Label and orient
its pseudolines arbitrarily thus getting a linking matrix $X$. Reverse the
orientations of all lines having an odd number of crossings with positive sign.
Call the resulting orientation {\it Eulerian}. It is unique up to 
reversing the orientation of all lines.

An {\it Eulerian linking matrix} $X_E$ associated to an Eulerian 
orientation of $D$ is uniquely defined up to conjugation by a 
permutation matrix. Its invariants
coincide with those of the switching class of $X_E$ but are slightly 
easier to compute since there are no sign ambiguities. In particular,
some of them can be computed using only $O(n^2)$ operations.

An {\it Eulerian partition} of the set of pseudolines of a diagram
consisting of an odd number of pseudolines is by definition
the partition of the pseudolines into subsets ${\mathcal L}_k$
consisting of all pseudolines involved in exactly $2k$ 
positive crossings for an Eulerian orientation.

\medskip

A few more invariants of Eulerian matrices are:

\begin{enumerate}

\item The total sum $\sum_{i,j} x_{i,j}$ of all entries in an Eulerian
linking matrix $X_E$ (this is of course equivalent to the computation of
the number of entries $1$ in $X_E$). The computation of this invariant
needs only $O(n^2)$ operations.

\item Its signature $\epsilon=\prod_{i<j} x_{i,j}$. The easy
identity
$$\epsilon=(-1)^{{n\choose 2}+\sum_{i<j}x_{i,j}}$$
relates the signature to the total sum $\sum_{i,j} x_{i,j}$ of all entries
in an Eulerian linking matrix.

\item The number of rows of $X_E$ with given row-sum (this
can also be computed using $O(n^2)$ operations). These numbers yield
of course the
cardinalities of the sets ${\mathcal L}_0,{\mathcal L}_1,\dots$.

\item All invariants of the associated
Eulerian graph (having edges associated to entries
$x_{i,j}=1$) defined by $X_E$, e.g. the number of triangles or of
other fixed subgraphs. 

\end{enumerate}

\medskip

For example, for $7$ vertices, there are $54$ different Eulerian
graphs, $36$ different sequences (up to a permutation of the vertices)
of vertex degrees, and $18$ different numbers for the cardinality of
1's in $X_E$.

\begin{figure}[h]
\epsfysize=4cm
\epsfbox{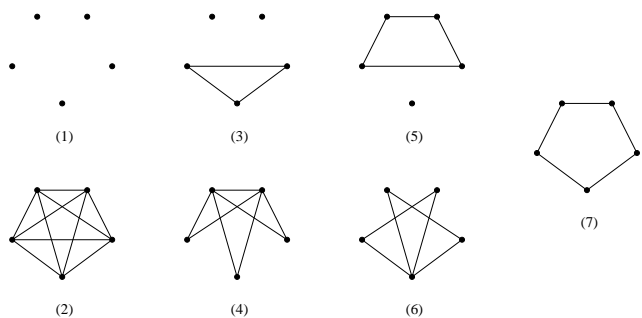}
\caption{All Eulerian graphs on $5$ vertices}\label{Euler_graphs_5}
\end{figure}

Figure \ref{Euler_graphs_5} shows all seven Eulerian graphs on $5$ vertices.

\subsection{Switching classes of even order - Euler partitions}\label{even_case}

The situation in this case is slightly more complicated and less 
satisfactory.

There exists a natural partition of the $2n$ rows $R$ of $X$ into two subsets
$R_+$ and $R_-$ according to the sign
$$\epsilon_i=\prod_{j\not=i}x_{i,j}\ \prod_{s<t}x_{s,t}$$
associated to the $i-$th row of $X$. This sign is indeed well-defined 
since switching (conjugation by a diagonal $\pm 1$ matrix $D$) 
multiplies both factors $\prod_{j\not=i}x_{i,j}$
and $\prod_{s<t}x_{s,t}$ by $\det (D)\in \{\pm 1\}$.

Since $\prod_i \epsilon_i=\prod_{i\not=j }x_{i,j}
\left(\prod_{s<t}x_{s,t} \right)^{2n}=1$, 
both subsets $R_+,R_-$ have
even cardinalities.

If $R_+$ (or equivalently, $R_-$) is non-empty, it defines a symmetric
submatrix $X_+$ of even size $\sharp(R_+)$ corresponding to all
rows and columns with indices in $R_+$. Iterating the above construction
we get thus a partition $R_+=R_{++}\cup R_{+-}$. This construction is
most conveniently encoded by a planar rooted binary tree which we
call the {\it Euler-tree} of $X$: Draw a root $R$ corresponding
to the row-set $R$ of $X$. If the partition $R=R_+\cup R_-$ is 
non-trivial join the root $R$ to  a left successor called 
$R_-$ and a right successor called $R_+$. The Euler tree of $X$ is now 
constructed recursively by gluing the root $R_\pm $ of the Euler tree 
associated to $X_\pm $ onto the corresponding successor $R_\pm $
of the root $R$.

The leaves of the Euler tree $T(X)$ of $X$ correspond to subsets 
$R_w$ (with $w\in\{\pm\}^*$) of even cardinality $2n_w$ 
summing up to $2n$. The leaves of $T(X)$ define
symmetric submatrices in $X$ which we call \lq\lq Eulerian'':
All their row-sums are identical modulo $2$ and can be chosen to be even
perhaps after switching an odd number of vertices (corresponding 
to conjugation by a diagonal $\pm 1$ matrix of determinant $-1$). 
The row partition $R=R_+\cup R_-$ of an Eulerian matrix $X$ 
of even size is by definition trivial. The sign $\epsilon\in\{\pm \}$
such that $R=R_\epsilon$ is called the the \lq\lq signature'' of $X$.
An Eulerian matrix of size $2$ has always signature $1$. For 
Eulerian matrices of size $2n\geq 4$ both signs can occur 
as signatures
since changing the signs of the entries $x_{i,j},\ 1\leq i\not= j\leq 3$
reverses the signature of an Eulerian matrix. 
The signature of an Eulerian matrix encodes the parity
of the number of edges in an Eulerian graph (having only vertices 
of even degrees) in the switching class of $X$.

\medskip

\noindent
{\bf Enumerative digression.} Associating a weight $n\in
\{1,2,\dots\}$ to a (non-empty) leaf corresponding to an Euler matrix
of order $2n$, the generating function $F(z)=
\sum_{n=0} \alpha_n z^n$ enumerating the number $\alpha_n$ of Euler
trees with total weight $n$ (associated to even switching classes
of order $2n$) satisfies the equation
$$F(z)=\frac{1}{1-z}+\left(F(z)-1\right)^2$$
(with $\alpha_0=1$ corresponding to the empty tree).
Indeed, Euler trees reduced to a leaf contribute $1/(1-z)$ to
$F(z)$. All other Euler trees are obtained by gluing two Euler
trees of strictly positive weights below a root and are enumerated
by the factor $\left(F(z)-1\right)^2$.

Solving for $F(z)$ we get the closed form
$$F(z)=\sum_{n=0}^\infty \alpha_nz^n=\frac{3(1-z)-\sqrt{(1-z)(1-5z)}}
{2(1-z)}.$$
showing that
$$\lim_{n\rightarrow\infty} \frac{\alpha_{n+1}}{\alpha_n}=5\ .$$
The first terms $\alpha_0,\alpha_1,\dots$ are given by
$$1,1,2,5,15,51,188,731,2950,12235,\dots$$
(see also Sequence A7317 in \cite{Sl}).

Similarly, the generating function $F_s(z)=
\sum_{n=0} \beta_n z^n$ enumerating the number $\beta_n$ of signed 
weighted Euler
trees (keeping track of the signature of all leaves with weight $\geq 2$)
with total weight $n$ satisfies the equation
$$F_s(z)=\frac{1+z^2}{1-z}+\left(F(z)-1\right)^2\ .$$
We get thus 
$$F_s(z)=\sum_{n=0}^\infty \beta_nz^n=\frac{3(1-z)-\sqrt{(1-z)(1-5z-4z^2)}}
{2(1-z)}.$$
and
$$\lim_{n\rightarrow\infty} \frac{\beta_{n+1}}{\beta_n}=\frac{5+\sqrt{41}}{2}
\sim 5.7016\ .$$
The first terms $\beta_0,\beta_1,\dots$ are given by
$$1,1,3,8,27,104,436,1930,8871,41916,\dots\ .$$

The leaves of the Euler tree define a natural partition of the set of
rows of $X$ into subsets. We call
this partition the {\it Euler partition}.
 
\begin{exa} \label{exampletree} The symmetric matrix
$$\left(\begin{array}{rrrrrrrrrr}
 0 & -1 &  1 & -1 &  1 &  1 &  1 &  1 &  1 & -1 \cr
-1 &  0 & -1 &  1 & -1 &  1 &  1 &  1 & -1 & -1 \cr
 1 & -1 &  0 & -1 &  1 &  1 &  1 & -1 & -1 &  1 \cr
-1 &  1 & -1 &  0 & -1 &  1 &  1 & -1 &  1 & -1 \cr
 1 & -1 &  1 & -1 &  0 &  1 & -1 & -1 &  1 &  1 \cr
 1 &  1 &  1 &  1 &  1 &  0 &  1 & -1 &  1 & -1 \cr
 1 &  1 &  1 &  1 & -1 &  1 &  0 & -1 & -1 &  1 \cr
 1 &  1 & -1 & -1 & -1 & -1 & -1 &  0 &  1 &  1 \cr
 1 & -1 & -1 &  1 &  1 &  1 & -1 &  1 &  0 &  1 \cr
-1 & -1 &  1 & -1 &  1 & -1 &  1 &  1 &  1 &  0
\end{array}\right)$$
yields the Euler partition
$$R_+=\{1,2,4,7,8,9\} \ \ \cup \ \ R_{-+}=\{6,10\}\ \ \cup\ \ R_{--}=\{3,5\}$$
(where the Eulerian submatrix associated to $R_+$ has signature $1$)
with Euler tree presented in Figure \ref{tree}.

\begin{figure}[h]
\epsfysize=3cm
\epsfbox{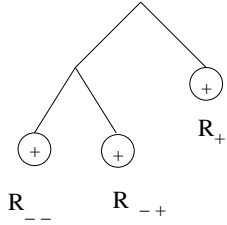}
\caption{The Euler tree of Example 
\ref{exampletree}}\label{tree}
\end{figure}

\end{exa}

\begin{exa}\label{exam_euler_even}
The two matrices mentioned in Remark \ref{character} are indeed not 
switching-equivalent: The row-partition $R=R_+\cup R_-$ of the 
first matrix is given by $R_+=\{1,3,4,6\}$ and $R_-=\{2,5,7,8\}$. 
The associated Euler
matrices $X_+$ and $X_-$ are both of signature $1$. On the other hand,
the second matrix is Eulerian and has also signature $1$.
\end{exa}

Let us mention a last invariant related to the
Euler tree for a switching class $X$ having even order $2n$.
Let $R_1,\dots,R_r\subset R$ be the Euler partition of $X$ . 
For $1\leq i\leq j\leq r$ define numbers $a_{i,j}
\in \{\pm 1\}$ by
$$a_{i,j}=\left\lbrace\begin{array}{ll}
\displaystyle 
\prod_{t\not= s_0\in R_i} x_{s_0,t}\prod_{s,t\in R_i,s<t} x_{s,t}
=\hbox{signature}(X_i)\quad&i=j\\
\displaystyle \prod_{s\in R_i,\ t \in R_j}x_{s,t}&i\not=j\end{array}
\right.
$$
where $s_0\in R_i$ is a fixed element (and where $X_i$ denotes
the Eulerian 
submatrix defined by $R_i$). 
One can easily check that the numbers $a_{i,j}$ are
well-defined.

This invariant has an even stronger analogue
for switching classes of odd order: Given an Eulerian matrix of order
$2n+1$ with Euler partition $A_0,\dots,A_r$ (where $A_i$
corresponds to the set of vertices of degree $2i$ in the Euler graph $\Gamma$)
one can consider the numbers 
$$a_{i,j}=\sum_{s\in A_i,t\in A_j} x_{s,t},\ 0\leq i,j$$
related to the number of edges joining vertices of given degree in $\Gamma$.

\section{Spindles}\label{spindles}

\subsection{Spindle-configuration}\label{sec_spindle_def}

Recall that 
a {\it spindle} is a configuration of skew lines 
intersecting an auxiliary line $A$, called its {\it axis}. 
A {\it spindle-configuration} 
(or a {\it spindle structure}) is a configuration of skew lines 
isotopic to a spindle.

The orientation of the axis $A$ induces a linear order $L_1<\dots<L_n$
on the $n$ lines of a spindle $C$. 
Each line $L_i\in C$ defines
a hyperplane $\Pi_i$ containing $L_i$ and the axis $A$. 

\begin{figure}[h]
\epsfysize=5cm
\epsfbox{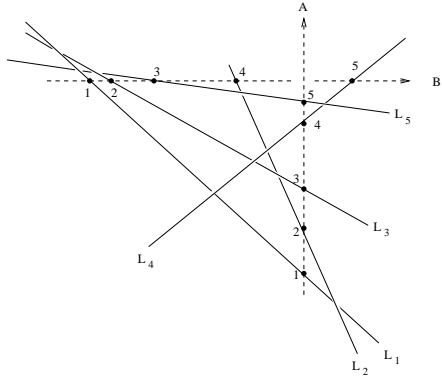}
\caption{A spindle}\label{spindle}
\end{figure}

A second directed 
auxiliary line $B$ (called a {\it directrix}) 
in general position with respect to $A,\Pi_1,\dots,\Pi_n$ and crossing 
$A$ negatively,
intersects the hyperplanes $\Pi_1,\dots,\Pi_n$ in points 
$\sigma(L_i)= B\cap \Pi_i$.
One can assume $\sigma(L_i)\in B$ by a suitable rotation of
$L_i\subset \Pi_i$ around $A\cap \Pi_i$.
Since the orientation of $B$ induces a linear order on the points 
$\sigma(L_i)$, we get a spindle-permutation (still denoted) $
i\longmapsto \sigma(i)$ of the set $\{1,\dots,n\}$ by identifying 
the two linearly ordered sets $L_1,\dots,L_n$ and $\sigma(L_1),\dots.
\sigma(L_i)$ in the obvious way with $\{1,\dots,n\}$. 
Figure \ref{spindle} displays an example corresponding  to
$\sigma(1)=1,\ \sigma(2)=4,\ \sigma(3)=2,\ \sigma(4)=5,\ \sigma(5)=3$.

Spindles can also be represented by (the isotopy classes of) 
configurations of skew lines with all lines contained in 
distinct affine horizontal planes of ${\mathbb R}^3$
(horizontal configurations). An associated spindle 
permutation encodes then the two orders on the set of lines 
given by considering the heights of the horizontal planes containing
them and by their slopes after vertical projection onto such a
horizontal plane.  

A linking matrix $X$ of a spindle $C$ is easily computed as follows.
Transform $C$ isotopically into a spindle with 
oriented axis $A$ and directrix $B$ as
above. Orient a line $L_i$ from $L_i\cap A$ to $\sigma(L_i)=L_i\cap B$. A 
straightforward
computation shows that the linking matrix $X$ of this labeled and 
oriented skew configuration has coefficients
$$x_{i,j}={\rm sign}((i-j)(\sigma(i)-\sigma(j)))$$
where $\hbox{sign}(0)=0$ and $\hbox{sign}(x)=\frac{x}{\vert x\vert}$
for $x\not=0$ and where
$\sigma$ is the corresponding spindle-permutation.

For example, the configuration of $5$ skew lines depicted in 
Figure \ref{spindle} corresponds to the linking matrix
$$X=\left(\begin{array}{rrrrr}
0 &  1 &  1 &  1 &  1\cr
1 &  0 & -1 &  1 & -1\cr
1 & -1 &  0 &  1 &  1\cr
1 &  1 &  1 &  0 & -1\cr
1 & -1 &  1 & -1 &  0\end{array}\right)\ .$$

\begin{rem} 
A configuration $C$ of $n$ skew lines 
has a spindle structure if and only if its
mirror configuration $\overline C$ has a spindle structure. A spindle 
permutation $\overline \sigma$ for $\overline C$
is then for instance given by
${\overline \sigma}(i)=n+1-\sigma(i),\ 1\leq i\leq n$, where 
$\sigma$ is a spindle-permutation for $C$.
\end{rem}
 
\subsection{Spindle-equivalent permutations}

The aim of this subsection is to show that three types of
 transformations of spindle-permutations, defined
by Crapo and Penne in \cite[Section 15]{CP}
preserve the associated spindle-configuration, up to isotopy. 
This is the easy direction in Theorem \ref{main} and follows
for instance from Theorem 3.2 in \cite{KM}.

\medskip

A (linear) {\it block} of size $k$ or a {\it $k-$block}  
in a permutation $\sigma$ of $\{1,\dots,n\}$ is a subset
$\{i+1,\dots,i+k\}$ of $k$ consecutive integers in $\{1,\dots,n\}$
such that 
$$\sigma(\{i+1,\dots,i+k\})=\{j+1,\dots,j+k\}$$
(i.e. the image under $\sigma$ of a set $\{i+1,\dots,i+k\}$ of
$k$ consecutive integers is again a set of $k$ consecutive integers).
In the sequel, we denote by $[\alpha,\beta]=\{\alpha,\alpha+1,\dots,
\beta-1,\beta\}\subset \{1,\dots,n\}$ a subset of consecutive integers and
by $\sigma([i+1,i+k])=[j+1,j+k]$ a $k-$block as above.

%

Recall that two spindle-permutations are {\it equivalent} 
(see \cite[Section 15]{CP})
if they are equivalent under the equivalence relation generated by
\begin{enumerate}
\item (Circular move) 
$$\sigma\sim \mu\hbox{ if }\mu(i)=(s+\sigma((i+t)\pmod n)) \pmod n$$ 
for some integers $0\leq s,t<n$ (all integers are modulo $n$).

\item (Vertical reflection of a block or local reversal) 
$\sigma\sim \mu$ if $\sigma([1,k])=[1,k]$ and
$$\mu(i)=\left\{\begin{array}{ll}
k+1-\sigma(k+1-i)\qquad &i\leq k\cr
\sigma(i)&i> k
\end{array}\right.$$
for some integer $k\leq n$ (see Figure \ref{vertical_ref}).

\begin{figure}[h]
\epsfysize=3cm
\epsfbox{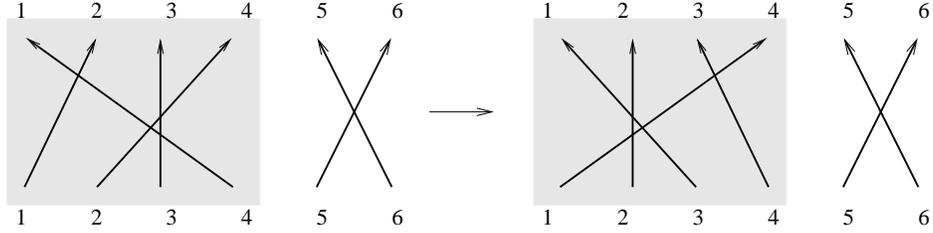}
\caption{Vertical reflection of a block}\label{vertical_ref}
\end{figure}                                                                 

\item  (Horizontal reflection of a block or local inversion) 
$\sigma\sim \mu$ if there exists an integer $1<k\leq n$ 
such that $\sigma([1,k])=[1,k]$ and
$$\mu(i)=\left\{\begin{array}{ll}
\sigma^{-1}(i)\qquad &i\leq k\cr
\sigma(i)&i> k\end{array}\right.$$
(see Figure \ref{horizontal_ref}).

\begin{figure}[h]
\epsfysize=3cm
\epsfbox{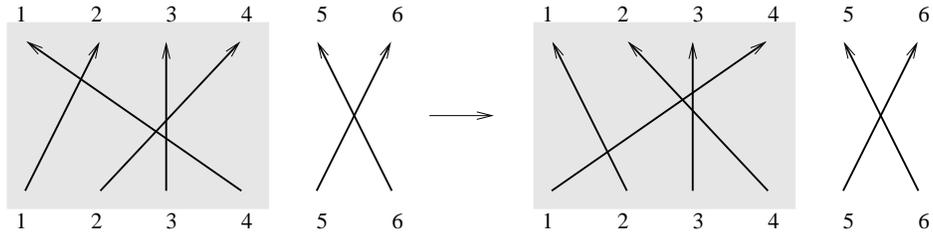}
\caption{Horizontal reflection of a block}\label{horizontal_ref}
\end{figure}                                                                 

\end{enumerate}

Circular moves suggest to extend the linear order on $\{
1,\dots,n\}$ to the cyclic order induced by the compactification 
$\R\P^2\supset \R$. We have an 
obvious notion of cyclic blocks and can consider linear horizontal
and vertical reflections of cyclic blocks related to
vertical and horizontal reflections as above by conjugations 
involving circular moves. Such (more general) moves lead
to the same equivalence relation as the 
three moves considered above and we allow them for the
sake of concision.

We restate and reprove the easy direction (see also \cite{KM})
of Theorem \ref{main}:

\begin{prop}\label{maineasy}
Equivalent spindle-permutations yield isotopic spindle-configurations.
\end{prop}

\begin{proof}
A transformation of type (1) amounts to pushing the last few lines
of the spindle on the axis and directrix through infinity. 
This can be done isotopically (continuously without leaving the 
set of skew configurations).

An isotopy inducing a vertical reflection (type (2) above)
can be described as follows:
Consider the two complementary subblocks
$\sigma([1, k]) = [1, k]$ and  
$\sigma([k+1, n]) = [k+1, n]$ in $\sigma$. All lines of the first 
block $\sigma([1, k]) = [1, k]$
can be squeezed isotopically into the interior of a small 
one-sheeted hyperboloid $H$ whose axis of revolution
intersects orthogonally the axis $A$ and the directrix $B$ of $\sigma$. 
Moreover, we may assume that no line of the complementary block
$\sigma([k+1, n]) = [k+1, n]$ intersects the interior of $H$.
The isotopy of skew lines given by rotating the interior 
(containing the block $\sigma([1, k]) = [1, k]$)
of $H$ by a half-turn around its axis of revolution
induces a vertical reflection (see Figure \ref{rule2}).

\begin{figure}[h]
\epsfysize=7cm
\epsfbox{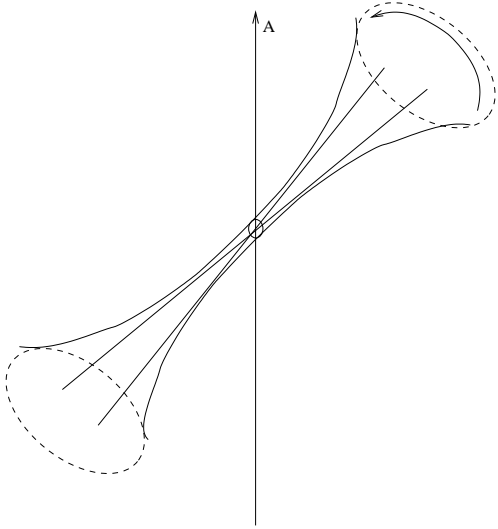}
\caption{Isotopy for a transformation of type (2)}\label{rule2}
\end{figure}                                                                 

For constructing an isotopy inducing a horizontal reflection (type (3)) 
we start as above by pushing
the $k$ lines of the first block $\sigma([1,k])=[1,k]$
into a small hyperboloid $H_1$ with
revolution axis $C_1$ intersecting the axis $A$ and the directrix $B$
orthogonally. Denote by $I_B$ the open
segment of $B$ contained in the interior of $H_1$.
Moreover, suppose that the directions of the axis $A$ 
and the directrix $B$ are orthogonal.
Push the lines of the complementary subblock $\sigma([k+1,n])=[k+1,n]$
in the positive sense along the directrix $B$ until they can be squeezed 
into the interior of a 
small revolution hyperboloid $H_2$ not intersecting
$H_1$ with revolution axis
$C_2$ parallel to the directrix $B$ of $\sigma$.

Rotate the interior of the first hyperboloid $H_1$ 
containing the subblock $\sigma([1,k])=[1,k]$
by a half-turn around the revolution axis $C_2$ of $H_2$
(see Figure \ref{rule3}) and call the resulting hyperboloid $H_1'$.
Finally, rotate the hyperboloid $H_1'$ and the lines inside it by 
$\pi/2$ around its revolution axis $C_1'$ and translate it along $C_1'$
until the image of $I_B$ is contained in 
the directrix $A$ of $\sigma$. 
This yields a spindle-configuration whose spindle-permutation
is related by a circular move and a horizontal reflection 
(and perhaps a vertical reflection of the first subblock, depending on
the sense of the half-turn around $C'_1$) 
to the initial spindle-permutation.

\begin{figure}[h]
\epsfysize=15cm
\epsfbox{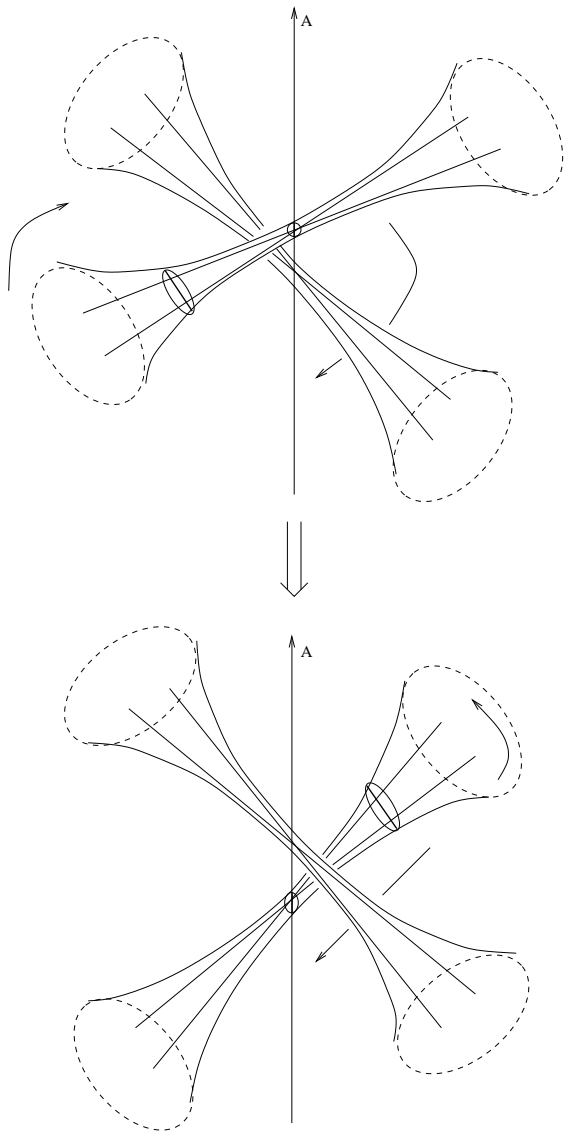}
\caption{Isotopy for realizing a transformation of type (3)}\label{rule3}
\end{figure}                                                                 

\end{proof}

\section{Isotopic spindle-configurations}\label{isospindles}

Our next result describes all blocks (up to spindle moves) in 
spindle-configurations.

\begin{prop} \label{Xsubblock} Let $\sigma$ be a 
spindle-permutation with linking matrix $X$.
Let $I$ be a subset of indices such that
$$x_{i_1,i_2}x_{i_1,j}x_{i_2,j}\in\{\pm 1\}$$ 
depends only on $i_1\not=i_2\in I$ and is independent of 
$j\not\in I$.

Then, up to spindle moves of $\sigma$, the lines $\{L_i\}_{i\in I}$ 
corresponding to $I$ form a block of $\sigma$.
\end{prop}

\begin{cor} \label{isosubblock} Let $B\subset \sigma$ be a block of a spindle-permutation $\sigma$. Let $\iota$ be an isotopy between $\sigma$ and 
a second spindle-permutation $\sigma'$. Then, up to spindle moves
of $\sigma'$, the set of lines $\iota(B)\subset \sigma'$ is a block.
\end{cor}

\begin{proof}[Proof of Proposition \ref{Xsubblock}] 
We suppose the rows and columns of 
$X$ indexed by $\{1,\dots,n\}$.

We consider first the case $\sharp(I)=2$. 
Up to circular moves, we can assume that $I=\{1,\alpha\}$ corresponds
to the lines $(1,1),(\alpha,\beta)\in\sigma$. Up to replacing 
$\sigma$ by its mirror spindle $\overline \sigma$, we can also assume
$x_{1,\alpha}x_{1,\gamma}x_{\alpha, \gamma}=1$ for 
$\gamma\in \{2,\dots,n\}\setminus\{\alpha\}$.
This implies $(\alpha-\gamma)(\beta-\sigma(\gamma))>0$
and shows $\alpha=\beta$.
Moreover, the lines $(2,\sigma(2)),(3,\sigma(3)),\dots,
(\alpha-1,\sigma(\alpha-1)),(\alpha,\alpha)$
form a block $B\subset \sigma$ and a vertical reflection with respect to $B$  
sends $(\alpha,\alpha)$ onto $(2,2)$.

Suppose now that the subset $I$ containing $k>2$ lines is
minimal in the sense that it contains no strict
subset $I'\subset I$ of $k'\geq 2$ indices satisfying the 
condition of Proposition \ref{Xsubblock}.
We claim that the corresponding subset of lines 
$\{(i,\sigma(i))\}_{i\in I}$ is a block of $\sigma$. 
Indeed, suppose that our claim does not hold.
Up to a horizontal reflection and a circular move, there exist
$1<j_1<i<j_2\leq n$ with $1,i\in I$ and $j_1,j_2\not\in I$. 
Up to replacing $\sigma$ by its mirror and up to circular moves
we can assume $1=\sigma(1)<\sigma(j_1)<\sigma(i)$. This implies
$\sigma(i)<\sigma(j_2)$ and shows that the two non-empty sets 
$I_1=I\cap\left(\{1,2,\dots, j_1-1\}\cup\{j_2+1,\dots,n\}\right)$ and 
$I_2=I\cap\{j_1+1,\dots,j_2-1\}$ yield a partition of $I$ into subsets
which satisfy both the condition of Proposition \ref{Xsubblock}. 
Since at least one
of the sets $I_1,I_2$ contains $\geq 2$ elements, we get a contradiction
with the assumed minimality of $I$.

In the general case we consider a subset $I'\subset I$ which satisfies the
condition of the Proposition \ref{Xsubblock} and which contains 
either only two lines or is minimal as defined previously. Since Proposition 
\ref{Xsubblock} holds for $I'$ we can suppose (up to spindle moves 
if $\sharp(I')=2$)  
that the lines of $I'$ form a block $B'\subset\sigma$. Considering the lines
of $B'$ as rigidly linked  (and thus allowing only spindle-moves of $\sigma$ 
transforming all lines of $B'$ similarly) 
we can consider $B'$ as being represented by a 
single line $L'\in B'$.
This yields a smaller spindle-permutation $\tilde \sigma$ 
and a subset of indices $\tilde I\subset $ 
satisfying the condition of Proposition \ref{Xsubblock} for the
corresponding linking matrix $\tilde X$ obtained by removing  from $X$ all
rows and columns corresponding to $B'\setminus L'$. 
Proposition \ref{Xsubblock} holds now for $\tilde \sigma$
by induction on the number of lines 
and gluing back the rigid block $B'$ (which is well-defined up to 
isotpy and a
vertical and horizontal reflection of $B'$) 
onto $L'\in\tilde\sigma$ yields the result.
\end{proof}

\begin{proof}[Proof of Corollary \ref{isosubblock}] The set of indices 
corresponding to a subblock $B\subset \sigma$ satisfies the condition of 
Proposition \ref{Xsubblock} for a common linking matrix $X$ of the
isotopic spindles $\sigma$ and $\sigma'$.
\end{proof}


\begin{thm} \label{mainiso} 
Let $\iota$ be an isotopy relating two spindle-configurations
(associated to spindle-permutations) $\sigma$ and $\sigma'$. Then the bijection
from the lines of $\sigma$ onto the lines of $\sigma'$ induced by
$\iota$ can be given by spindle-moves. \end{thm}

\begin{proof}[Proof of Theorem \ref{main}] Follows from Proposition 
\ref{maineasy} and Theorem \ref{mainiso}.
\end{proof}

\begin{rem} There might
exist \lq\lq exotic'' isotopies between spindle-equivalent 
spindle-configurations which do not arise
from (a continuous deformation of a sequence of) spindle-moves. 
\end{rem}
 
A few notations: We generalize the notions of spindle-permutations,
spindle-configurations etc. as follows:
A {\it spindle-permutation} is a bijection $\sigma: E\longrightarrow F$
between two finite subsets $E,F\subset {\mathbb R}$ which we consider
either linearly ordered or cyclically ordered by the cyclic order
induced on ${\mathbb R}$ from the compactification
${\mathbb R}\subset {\mathbb R}\cup\{\infty\}={\mathbb R}P^1\sim 
{\mathbf S}^1$. For $e\in E$ we denote by $(e,\sigma(e))$ the line
of the spindle-configuration associated in the obvious way to $\sigma$
and we identify $\sigma$ with the set $\{(e,\sigma(e)\}_{e\in E}$ of
its lines. For $e\in E$, the 
notation $\sigma\setminus (e,\sigma(e))$ denotes
the spindle or spindle-permutation obtained by restricting $\sigma$ to 
$E\setminus \{e\}$.
For $e_1<e_2\in E$ and $f_1<f_2\in F$ we denote by $[e_1,e_2]$
the subsets $[e_1,e_2]\cap E$ and $[f_1,f_2]\cap F$ of $E$ and $F$.
For subsets $E'\subset E,F'\subset F$ of the same cardinality
such that $\sigma(e)\in F'$ for all $e\in E'$ we 
denote by $(E',\sigma(E')=F')$ the spindle-permutation obtained 
by restricting $\sigma$ to $E'$. A subblock of $\sigma$ 
can thus be written as
$([e_1,e_2],\sigma([e_1,e_2])=[f_1,f_2])\subset \sigma$ and we will also 
use the shorthand notation $\sigma([e_1,e_2])=[f_1,f_2]$. In the sequel,
a $k-$block (of a spindle-permutation $\sigma$) will almost always 
denote a cyclic block consisting of $k$ lines, 
i.e. a subset $E'\subset E$ of $k$ cyclically consecutive elements 
with $\sigma(E')$ cyclically consecutive in $F$.
Let us also remark that given a spindle-permutation $\sigma:E\longrightarrow
F$, its mirror configuration is for instance associated to 
the spindle-permutation $\overline \sigma:E\longrightarrow \overline F$
where $\overline F =F$ as a set but equipped with the opposite (cyclic) order.
The application $\sigma\longmapsto \overline \sigma$ which replaces
a spindle-permutation by its mirror enjoys good properties (preserves
the spindle-equivalence relation, the isotopy relation, 
yields a bijection between subblocks etc.) and will
often be used to reduce the number of possible cases.
 
The proof of Theorem \ref{mainiso} is by induction 
on the number $n$ of lines involved in $\sigma$ and $\sigma'$.
The result holds clearly for configurations of $\leq 3$ lines
(in this case, spindle-moves generate the complete permutation group
of all lines).

Call a permutation {\it irreducible} if it contains no
non-trivial block (consisting of $2\leq k\leq n-2$ 
cyclically consecutive lines).

Call a block {\it minimal} if it consists of $2\leq k\leq n-2$ lines
and if it contains no subblock of strictly smaller cardinality $k'\geq 2$.

Proposition \ref{Xsubblock} shows that the set of possible subblocks 
of a spindle-permutation $\sigma$ is encoded in its linking matrix.
Thus, either both or none of the spindle
permutations $\sigma,\sigma'$ are irreducible.

The proof of Theorem \ref{mainiso} splits into two cases,
depending on the reducibility of $\sigma$ and $\sigma'$.

\subsection{The reducible case}\label{reducible_case}

Consider a non-trivial $k-$subblock
$B\subset \sigma$. Corollary \ref{isosubblock} shows that,
up to spindle moves, $B'=\iota(B)$ is a non-trivial $k-$subblock of $\sigma'$.

Up to circular spindle moves we can assume $B$ and $B'$ be given by 
$\sigma([1,k])=[1,k]$ and $\sigma'([1,k])=[1,k]$. We denote by $\overline B=
\sigma\setminus B$ and $\overline B'=\sigma'\setminus B'$ 
the complementary blocks.

By induction on the number of lines, the bijection of lines obtained
by restricting $\iota$ to the subspindles $(1,\sigma(1))\cup\overline B$ 
and  $(\iota(1),\sigma'(\iota(1)))\cup \overline B'$ can be obtained
by an isotopy $\mu$ given by a composition of
spindle-moves. Up to a vertical and/or horizontal
reflection of $B$, the isotopy $\mu$
can be extended uniquely to all lines of 
$\sigma$ by considering the subblock $B\subset \sigma$ as rigid (and thus by 
allowing only spindle moves having the same effect on all lines of $B$).
Replacing $\iota$ with the isotopy $\mu^{-1}\circ \iota$ we can thus
suppose that the permutation induced by 
$\iota$ fixes all lines of $\overline{B}=\sigma\setminus B$. 
Applying the above
argument to the complementary block $\overline B$ we get the
result. \hfill $\Box$


\subsection{The irreducible case}\label{irreducible_case}

The case where $\sigma$ and $\sigma'$ are irreducible is more involved.
It splits into the three following subcases:

Subcase (1): There exists a line $L\in \sigma$ with $\sigma\setminus
L$ irreducible.

Subcase (2): There exists a line $L$ such that $\sigma\setminus L$ 
contains a minimal block $B$ consisting of $3\leq k<n-3$ lines.

Subcase (3): $\sigma\setminus L$ contains a $2-$block for every line 
$L\in \sigma$.

Subcases (1) and (2) are dealt with by induction on the number of lines.
We call a spindle-permutation giving rise to subcase (3) {\it exceptional}.
Subcase (3) is then handled by classifying all exceptional 
spindle-permutations.

\subsubsection{Subcase (1)}
Choose a line $L\in \sigma$ with $\sigma\setminus L$ irreducible.
By Corollary \ref{isosubblock}, the spindle-permutation
$\sigma'\setminus \iota(L)$ is also  irreducible.
Up to circular moves of $\sigma$ and $\sigma'$ we can
assume $L=(n,n)\in\sigma , \iota(L)=(n,n)\in\sigma'$.

By induction on $n$, the restriction of $\iota$ to the spindle
configurations $\sigma([1,n-1])=[1,n-1],\ \sigma'([1,n-1])=[1,n-1]$
yields a bijection between their lines 
which can be realized by spindle moves. 
Up to applying horizontal and vertical reflections to $\sigma'$,
there exist (by irreducibility of $\sigma\setminus L$ and 
$\sigma'\iota(L)$)
integers $0\leq \alpha,\beta<n-1$ such that
$$\sigma'(i)=\sigma\Big(i-\alpha\pmod {(n-1)}\Big)+\beta\pmod {(n-1)}$$
for $i=1,\dots,n-1$
where $x\pmod {(n-1)}\in\{1,\dots,n-1\}$. Say that line 
$(i,\sigma(i))$ (with $i<n$)
is {\it moved through infinity} if either $i\geq n-\alpha$ or
$\sigma(i)\geq n-\beta$. If both inequalities hold, we say that
$(i,\sigma(i))$ is {\it moved twice through infinity}. 

If there exists a line $(
i,\sigma(i))$
which is moved exactly once through infinity, then every line
$(j,\sigma(j)),\ 1\leq j\leq n-1$ is moved exactly once through infinity: 
Otherwise , consider a line $(i,\sigma(i))$ which is moved once and a
line $(j,\sigma(j))$ which is not moved through infinity or moved twice. 
This implies
that the $3-$subspindles  
$$\{(i,\sigma(i)),(j,\sigma(j)),(n,n)\}\hbox{  and } 
\{\iota(i,\sigma(i)),\iota(j,\sigma(j)),
\iota(n,n)\}$$ are mirrors (and thus not isotopic) which is 
impossible. Every line of $\sigma\setminus
(n,n)$ is thus moved through infinity exactly once which implies $\alpha+\beta=
n-1$. If $2\leq \alpha\leq n-3$, we get a contradiction with irreducibility of $\sigma\setminus(n,n)$. If $\alpha\in\{1,n-2\}$, we get a contradiction with irreducibility of $\sigma$. 

We can now assume that every line of $\sigma\setminus(n,n)$ is moved an
even number of times through infinity. This implies $\alpha=\beta$ and the existence  of a non-trivial subblock $\sigma([n-\alpha,n-1])=[n-\alpha,n-1]$
if $2\leq \alpha\leq n-3$ which contradicts irreducibility of $\sigma\setminus
(n,n)$. The case $\alpha=\beta\in\{1,n-2\}$ leads to a contradiction with irreducibility of $\sigma$. We get thus $\alpha=\beta=0$ which shows 
$\sigma=\sigma'$ and proves the result.

\subsubsection{Subcase (2)} 
Up to circular moves, we can assume that $\sigma\setminus (n,n)$ 
contains a block $B$ of size $3\leq k\leq n-4$.
Irreducibility of $\sigma$ shows now the existence
(up to horizontal and vertical reflections of $\sigma$)
of integers $\alpha,\beta\geq 1$ 
with $3\leq \alpha+\beta=\sharp(B)\leq n-4$ and of an integer $\gamma$
with $1\leq \gamma<n-2-\alpha-\beta$ such that
$$
B:\sigma([1,\alpha])\cup[n-\beta,n-1])= [\gamma+1,\gamma+\alpha+\beta]\ .$$
Up to considering mirror-configurations, we can suppose 
$\gamma<\sigma(\alpha)<\sigma(n-\beta)\leq \alpha+\beta+\gamma$.
There exist now $\alpha<\nu_1<\nu_2<n-\beta$ such that 
$\sigma(\nu_1)>\alpha+\beta+\gamma$ and $\sigma(\nu_2)\leq \gamma$. Indeed, 
otherwise $\sigma([\alpha+1,\alpha+\gamma])= [1,\gamma]$ and
$\sigma([\alpha+\gamma+1,n-\beta-1])=[\alpha+\beta+\gamma+1,n-1]$
which contradicts the irreducibility of $\sigma$ since this yields at
least one non-trivial block in $\sigma$.

\begin{figure}[h]
\epsfysize=3.5cm
\epsfbox{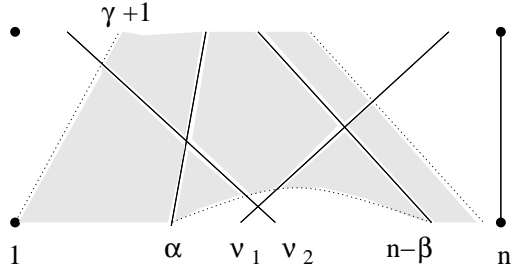}
\caption{A schematical picture of the
 subspindle $\tilde \sigma\subset \sigma$.}\label{pantalon}
\end{figure}                                                              

\begin{lem} The subspindle $\tilde \sigma\subset \sigma$ defined by
$$B\cup (\nu_1,\sigma(\nu_1))
\cup(\nu_2,\sigma(\nu_2))\cup (n,n)$$
is irreducible. 
\end{lem}

\begin{proof} Consider the $5-$subspindle $\tau\subset \tilde \sigma$ 
containing the five lines
$$(\alpha,\sigma(\alpha)),(\nu_1,\sigma(\nu_1)),(\nu_2,\sigma(\nu_2)),
(n-\beta,\sigma(n-\beta)),(n,n)\ .$$
The inequalities $\alpha<\nu_1<\nu_2<n-\beta<n$ and $\sigma(\nu_2)<\sigma(\alpha)<\sigma(n-\beta)<\sigma(\nu_1)<n$ imply easily that $\tau$ is
irreducible (it is enough to check that $\tau$ contains no $2-$block).
Thus any block $\tilde B\subset\tilde \sigma$ 
intersects $\tau$ in
a subset with cardinality $\sharp(\tilde B\cap \tau)\in
\{0,1,4,5\}$. Up to replacing 
$\tilde B$ by its complementary block $\tilde \sigma\setminus \tilde B
\subset \tilde \sigma$, we can assume $\sharp(\tilde B\cap \tau)\leq 1$.

If $\tilde B\cap \tau\subset \{(\alpha,\sigma(\alpha)),\ 
(n-\beta,\sigma(n-\beta))\}$ then $\tilde B\subset B$ is 
also a nontrivial subblock in $\sigma$. This contradicts the
irreducibility of $\sigma$.

If $\tilde B\cap \tau=\{(n,n)\}$, then a non-trivial subblock
$\tilde B$ of of $\tilde \sigma$ contains also at least one line of the set 
$\{(\nu_1,\sigma(\nu_1)),(\nu_2,\sigma(\nu_2))\}$. This contradicts
our assumption $\sharp(\tilde B\cap \tau)\leq 1$.

If $\tilde B\cap \tau=\{(\nu_1,\sigma(\nu_1))\}$, then a
non-trivial subblock $\tilde B$ contains also at least 
one line of $\{(\alpha,\sigma(\alpha)),(\nu_2,\sigma(\nu_2))\}$
contradicting again 
$\sharp(\tilde B\cap \tau)\leq 1$.

The case $\tilde B\cap \tau=\{(\nu_2,\sigma(\nu_2))\}$ is
analogous.

The case $\tilde B\cap \tau=\emptyset$ implies that $\tilde B$ is a
block of $\sigma$ which is impossible.
\end{proof}

Apply now Theorem \ref{mainiso} to $\iota$ restricted to
the subspindles $\tilde\sigma$ and $\tilde \sigma'=
\iota(\tilde\sigma)\subset \sigma'$ (which contain at most $n-1$ lines).
Assuming minimality of $B$, the subset $B'=\iota(B)$ is already a 
subblock in $\sigma'\setminus \iota(n,n)$ 
(see the proof of Proposition \ref{Xsubblock}).
Moreover, up to spindle moves, the relative position of
the subblock $B'\subset \sigma'\setminus (n,n)$ 
(with $(n,n)=\iota(n,n)$) inside $\sigma'$
is described by the integers 
$\alpha,\beta,\gamma$ considered above: One can indeed compute 
these integers by counting isotopy classes of suitable triplets $L_B\in
B, L_{\overline B}\in \sigma\setminus (B\cup (n,n)), (n,n)$ which are 
in bijection with the corresponding triplets in $\sigma'=\iota(\sigma)$. 
Induction on $n$
implies now, that (perhaps up to a vertical reflection of $\sigma'$) 
the permutation induced by $\iota$
fixes the lines of $B\cup (n,n)$ which can
be assumed to be common to $\sigma$ and $\sigma'$. 
Replacing $B$ by its complement $\overline B$
in $\sigma\setminus (n,n)$ we can find an irreducible subspindle 
$\overline \sigma\subset \sigma$ containing $\overline B, (n,n)$ and two 
suitable lines of $B$. Theorem \ref{mainiso} holds by induction on $n$
for the
isotopy $\overline \iota$ obtained by restricting $\iota$ to the subspindles
$\overline \sigma\subset \sigma,\ \iota(\overline \sigma)\subset \sigma'$ 
and implies that $\iota$ fixes also all lines of $\overline B$. (In this case, 
we have no longer to care about minimality of $\overline B$: the 
corresponding parameters $\overline \alpha,\overline \beta,\overline \gamma$ 
are fixed by the relative 
position of the already coinciding subblocks $B=B'$ 
inside $\sigma$ and $\sigma'$.)

This shows $\sigma=\sigma'$ and ends the proof of subcase (2).

\subsubsection{Exceptional irreducible spindles}
Call a spindle-permutation $\sigma$
of $n\geq 4$ lines {\it exceptional} 
if $\sigma$ is irreducible and
$\sigma\setminus L$ contains a (cyclic) $2-$block for every line 
$L\in \sigma$.

\begin{prop} (i) For $n\geq 5$ odd, the spindle-permutation 
$\tau=\tau_n$ of $\{0,\dots,n-1\}$ defined by 
$$\tau:i\longmapsto \tau(i)=2i\pmod n,\ 0\leq i\leq n-1$$
and its mirror $\overline \tau=\overline \tau_n$ given (up to a circular move)
by 
$$\overline \tau:i\longmapsto \tau(i)=-2i\pmod n,\ 0\leq i\leq n-1$$
are exceptional. The spindle-permutations $\tau_5$ and $\overline \tau_5$
are spindle-equi\-valent. For $n>5$ odd, the spindle-permutations
$\tau_n$ and $\overline \tau_n$  are not spindle-equivalent and have
linking matrices which are not in the same switching class. In particular,
the associated spindle-configurations are non-isotopic.

\ \ (ii) The spindle-permutations $\tau_n,\ \overline\tau_n,\ n\geq 5$ odd,
are the only exceptional spindle-permutations having $\geq 4$ lines, 
up to spindle-equivalence.

\ \ (iii) If $\iota$ is an isotopy of the exceptional
spindle-configuration $\tau_n$ onto
itself, then the line permutation $L\longmapsto \iota(L)$ induced by
$\iota$ can be realized by spindle moves.
\end{prop}

\begin{proof} We write $\tau=\tau_n$ for $n\geq 5$ odd.
We have $\tau=\tilde \tau^k$
with 
$$\tilde \tau^k:i\longmapsto \tilde\tau^k(i)=
2k+\tau(i-k\pmod n)\pmod n$$ showing that $\tau$ has a group of
automorphisms acting transitively on its lines. Since 
$\tau\setminus (0,0)$ contains the (cyclic) $2-$block $\{(\frac{n-1}{2},n-1),
(\frac{n+1}{2},1)\}$, we get thus a cyclic $2-$block in $\tau\setminus 
L$ for any line $L\in\tau$.

Suppose now that $\tau$ is reducible and consider a non-trivial subblock 
$B\subset \tau$. Up to replacing $B$ by the complementary block 
$\tau\setminus B$, 
we can assume that $B$ contains fewer than $n/2$ lines. Up to cyclic 
moves we can assume that $(0,0),(1,2)\in B$. This implies either 
$B=\tau\setminus (\frac{n+1}{2},1)$ which contradicts non-triviality of $B$
or $B$ contains the line $(\frac{n+1}{2},1)$. But then $B$ contains
either all lines $(k,2k),\ 1\leq k\leq \frac{n-1}{2}$ or all lines
$(k,2k-n),\ \frac{n+1}{2}\leq k\leq n-1$. Since $(0,0)\in B$,
we have in both cases $\sharp(B)\geq \frac{n+1}{2}$ which contradicts the
assumption $\sharp (B)<\frac{n}{2}$.

This shows that $\tau=\tau_n$ and its mirror $\overline \tau_n$ are 
exceptional. 

A vertical reflection transforms $\tau_5$ into $\overline\tau_5$
which proves their equivalence under spindle moves. 
For $n>5$ odd, $\tau_n\setminus L$ and
$\overline \tau_n\setminus \overline L$ contain both a unique $2-$block
$B$, resp. $\overline B$ 
(the choice of the lines $L,\overline L$ is irrelevant
since they are transitively permuted by automorphisms). We
get thus $3-$subspindles $L\cup B\subset \tau_n$ and $\overline L\cup 
\overline B\subset \overline \tau_n$ which are not isomorphic
since each is the mirror of the other. 
Since such a $3-$spindle and its isotopy 
class can be recovered from the linking matrix, the associated crossing 
matrices are not switching equivalent and the corresponding spindle
configurations are non-isotopic. This
proves assertion (i). 

Notice that this argument fails for $n=5$: In this case,
$\tau\setminus L$ gives rise to two complementary $2-$blocks 
$B,\overline B$ such that
$L\cup B$ and $L\cup\overline B$ are non-isotopic $3-$spindles.

An inspection shows that no irreducible $4-$spindle exists. This 
proves assertion (ii) for $n=4$.
Hence, we can suppose $n\geq 5$. Indexing the $n$ lines of an
exceptional spindle-permutation $\sigma$ by
$\{0,1,\dots, n-1\}$ and using circular moves, we can assume 
$\sigma(0)=0$. Denote by $B_0$ a $2-$block contained in
$\sigma\setminus (0,0)$. Up to a horizontal reflection, there exists
$1\leq \alpha\leq n-1$ such that $B_0=\{(\alpha,\sigma(\alpha)),
(\alpha+1,\sigma(\alpha+1))\}$ with
$\{\sigma(\alpha),\sigma(\alpha+1)\}=\{1,n-1\}$. Up to replacing $\sigma$
by its mirror spindle, we can suppose $\sigma(\alpha)=
n-1$ and $\sigma(\alpha+1)=1$. Irreducibility of $\sigma$ and 
$n\geq 5$ imply then $1<\alpha<n-2$.

For a $2-$block $B_{\alpha+1}$ contained in $\sigma\setminus(\alpha+1,1)$,
we have now the following three possibilities:
$$B_{\alpha+1}=\left\lbrace\begin{array}{ll}
\displaystyle \{(\alpha,n-1),(\alpha+2,n-2)\}&\hbox{case (a)}\\
\displaystyle \{(0,0),(n-1,2)\}&\hbox{case (b)}\\
\displaystyle \{(0,0),(1,2)\}&\hbox{case (c)}\end{array}\right.$$
Indeed, such a $2-$block is either given by the two lines
$(\alpha,x),(\alpha+2,y)$ (case (a)) 
with $x,y$ cyclically adjacent, or it is given by $(x,0),(y,2)$ 
(cases (b) and (c)) with $x,y$ cyclically adjacent, see Figure
\ref{abc} (with dots indicating the extremities and fatter lines
indicating the $2-$block $B_{\alpha+1}$).

\begin{figure}[h]
\epsfysize=2.5cm
\epsfbox{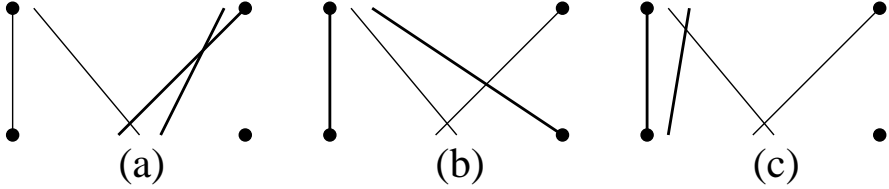}
\caption{The subcases (a), (b) and (c).}\label{abc}
\end{figure}

The following chart clarifies the steps of the proof:
$$\left\lbrace\begin{array}{ll}
a&\left\lbrace\begin{array}{ll}
aa\quad &\tau_5\hbox{ or contradiction}\\
ab\quad &\tau_5\hbox{ or contradiction}\end{array}\right.\\
b&\left\lbrace\begin{array}{ll}
ba\quad &\hbox{contradiction}\\
bb\quad &\left\lbrace\begin{array}{l}\hbox{contradiction}\\
\dots\hbox{ contradiction}\end{array}\right.\end{array}\right.\\
c&\left\lbrace\begin{array}{ll}
ca\quad &\tau_5\hbox{ or contradiction}\\
cb\quad &\left\lbrace\begin{array}{l}\tau_5\\
\left\lbrace\begin{array}{l}\tau_7\\\ddots\end{array}\right.\end{array}
\right.\end{array}\right.\end{array}\right.$$

Case (a): Consider a $2-$subblock $B_\alpha$ contained in $\sigma\setminus (\alpha,n-1)$. There are two possibilities for $B_\alpha$:
$$B_\alpha=\left\lbrace\begin{array}{ll}
\displaystyle \{(0,0),(\alpha+2,n-2)\}&\hbox{subcase (aa)}\\ 
\displaystyle \{(\alpha-1,2),(\alpha+1,1)\}&\hbox{subcase (ab)}\end{array}\right.$$ 

Subcase (aa) leads to $\alpha=n-3$
and the lines $\{(0,0),(n-3,n-1),(n-2,1),(n-1,n-2)\}$ form a non-trivial
subblock of $\sigma$ for $n>5$ which is a contradiction. For $n=5$ we get 
$\sigma=\tau_5$. 

Subcase (ab): The $2-$block $B_{\alpha-1}$ contained in 
$\sigma\setminus (\alpha-1,2)$ is either given by $B_{\alpha-1}=\{(\alpha+1,1),(\alpha+2,n-2)\}$ which implies $1+2=3=n-2$ and $\sigma=\tau_5$
or we have 
$B_{\alpha-1}=\{(\alpha,n-1),(0,0)\}$ which implies $\alpha=2$.
For $n=5$ we get again $\sigma=\tau_5$ and for $n>5$ the spindle
permutation $\sigma$ contains the non-trivial subblock
$\{(0,0),(1,2),(2,n-1),(3,1)\}$. This contradicts the irreducibility of
$\sigma$ and finishes the discussion of case (a).

Case (b): For the subblock $B_\alpha$ contained in
$\sigma\setminus (\alpha,n-1)$ we have
$$B_\alpha=\left\lbrace\begin{array}{ll}
\displaystyle \{(0,0),(\alpha+1,1)\}&\hbox{subcase (ba)}\\ 
\displaystyle \{(0,0),(1,n-2)\}&\hbox{subcase (bb)}\end{array}\right.$$ 

Subcase (ba) implies $\alpha=1$ which contradicts the irreducibility 
of $\sigma$.

Subcase (bb): The $2-$subblock $B_1$ associated to 
$\sigma\setminus
(1,n-2)$ is either given by $B_1=\{(0,0),(\alpha,n-1)\}$ which implies 
$\alpha=2$ contradicts the irreducibility (consider the $2-$block 
$\{(1,n-2),(2,n-1)\}\subset\sigma$) or $B_1$ is given by
$B_1=\{(\alpha,n-1),(\alpha-1,n-3)\}$. Considering the $2-$block of
$\sigma\setminus (\alpha-1,n-3)$ we get the existence of the line $(2,n-4)$
and iterating this argument, we get the existence of the lines 
$(i,n-2i),(\alpha+1-i,n+1-2i),\ 1\leq i\leq \lfloor\frac{\alpha+1}{2}\rfloor$
forming a non-trivial subblock
$([1,\alpha],\sigma([1,\alpha])=[n-\alpha,n-1])\subset \sigma$ 
and contradicting the irreducibility of $\sigma$.
This rules out case (b).

Case (c): For the $2-$block $B_\alpha$ of 
$\sigma\setminus (\alpha,n-1)$ we have
$$B_\alpha=\left\lbrace\begin{array}{ll}
\displaystyle \{(1,2),(\alpha+1,1)\}&\hbox{subcase (ca)}\\ 
\displaystyle \{(0,0),(n-1,n-2)\}&\hbox{subcase (cb)}\end{array}\right.$$ 

Subcase (ca) yields $\alpha=2$ and the four lines
$(0,0),(1,2),(2,n-1),(3,1)$ form a non-trivial subblock 
of $\sigma$ if $n>5$ which is a contradiction. 
For $n=5$ we get $\sigma=\tau_5$.

Subcase (cb): Considering $\sigma\setminus(1,2)$ we have either $n=5,\ \alpha=2$ and $\sigma=\tau_5$ or we obtain the existence of the line
$(\alpha+2,3)$. A symmetric argument (involving $\sigma\setminus (n-1,n-2)$)
yields the line $(\alpha-1,n-3)$.
We have now either $n=7,\ \alpha=3$ and $\sigma=\tau_7$ or we get the 
existence of two new lines $(2,4)$ and 
$(n-2,n-4)$ by considering the $2-$blocks of $\sigma\setminus (\alpha+2,2)$,
respectively $\sigma\setminus (\alpha-1,n-2)$. More generally, we get by
iteration of this construction either a $(4k+1)-$subspindle 
$\tilde \sigma_{4k+1}$ containing the lines 
$$(\alpha+1-k,n+1-2k),(\alpha+k,2k-1),(k,2k),(n-k,n-2k)$$
or a 
$(4k-1)-$subspindle 
$\tilde \sigma_{4k-1}$ containing the lines 
$$(k-1,2k-2),(n+1-k,n+2-2k),(\alpha+1-k,n+1-2k),(\alpha+k,2k-1)\ .$$

First consider $\tilde\sigma_{4k+1}$. If $\tilde\sigma_{4k+1}=\sigma$,
we get $\alpha=2k$ implying $\sigma=\tau_{4k+1}$.
Otherwise, a consideration of the $2-$blocks contained in $\sigma\setminus 
(k,2k)$ and $\sigma\setminus (n-k,n-2k)$ implies the existence 
of the lines $(\alpha+1+k,2k+1),(\alpha-k,n-1-2k)$ 
showing that $\tilde \sigma_{4k+1}$ can be 
extended to $\tilde \sigma_{4(k+1)-1}$.

Now consider $\tilde \sigma_{4k-1}$. If $\tilde \sigma_{4k-1}=\sigma$ we have
$\alpha=2k-1$ and $\sigma=\tau_{4k-1}$. Otherwise, considering the
$2-$blocks contained in $\sigma\setminus(\alpha+1-k,n+1-2k)$ 
and $\sigma\setminus(\alpha+k,2k-1)$, we get the existence of the lines 
$(n-k,n-2k)$ and $(k,2k)$. This shows that $\tilde \sigma_{4k-1}$ can 
be extended to $\tilde \sigma_{4k+1}$. Iteration of this construction
stops by finiteness of $\sigma$ and ends always with $\sigma=\tau_{n}$
for odd $n$. This proves assertion (ii).

Up to considering mirror configurations, it is enough to prove assertion
(iii) for $\tau_n$.
Up to spindle-moves (circular moves, horizontal and vertical reflections)
we can assume that the line-bijection induced by $\iota$ fixes the three
lines $(0,0),(\frac{n+1}{2},1),(\frac{n-1}{2},n-1)$. This implies that
it has also to fix the lines $(1,2)$ and $(n-1,n-2)$. Iteration
shows finally that the lines of all subspindles $\tilde \sigma_{4k+1}$
and $\tilde \sigma_{4k-1}$ are fixed. This implies the result.
\end{proof}

\begin{proof}[Proof of Theorem 1.3] Follows from subsections 
\ref{reducible_case}
and \ref{irreducible_case}
\end{proof}

\section{Theorem 3.2 of Kashin-Mazurovskii}\label{KaMa}

Theorem 3.2 of \cite{KM} can be restated in our terminology as
follows:
 
\begin{thm}\label{KM3.2} Spindle-configurations with switching-equivalent
linking-matrices are isotopic.
\end{thm} 

The proof is constructive and provides an explicit isotopy between 
the two spindle configurations.

Let $\sigma$ and $\mu$ be two spindles with switching-equivalent linking
matrices $X_\mu=P^tX_\sigma P$
where $P$ is a signed permutation matrix inducing
a bijection between the lines of $\sigma$ and $\mu$. 
(Recall that the linking matrix $X_\tau$ of a spindle with
permutation $\tau$ is defined by
$$(X_\tau)_{i,j}=\hbox{sign}\big((i-j)(\tau(i)-\tau(j))\big)\ \hbox{.)}$$
Up to circular moves,
we can assume that the line-bijection 
$$(i,\sigma(i))\longmapsto (i',\mu(i'))$$
induced by $P$ sends the first 
line $(1,\sigma(1)=1)$
of $\sigma$ onto the first line $(1,\mu(1)=1)$ of $\mu$. This implies 
that the conjugating matrix $P$ is an ordinary permutation matrix
and we have thus
$$(i-j)(\sigma(i)-\sigma(j))(i'-j')(\mu(i')-\mu(j'))>0$$
for all $1\leq i<j\leq n$.

For $t\in \R$ fixed, consider the line $L_i(t)=sP_\omega^i(t)+(1-s)P_\alpha^i(t)$ parametrized by $s\in \R$, oriented from 
$P_\alpha^i(t)$ to $P_\omega^i(t)$ where
$$\left\lbrace\begin{array}{lcl}
\displaystyle P_\alpha^i(t)&=&(1-t)(i,1,0)+t(0,1,i')\\
\displaystyle P_\omega^i(t)&=&(1-t)(0,-1,\sigma(i))+t(-\mu(i'),-1,0)
\end{array}\right.$$

\begin{prop} For $t\in\R$, the lines ${\mathcal L}(t)=\{L_1(t),\dots,L_n(t)\}$
form a skew configuration.

Moreover, ${\mathcal L}(0)$ realizes the spindle $\sigma$ and  
${\mathcal L}(1)$ realizes the spindle $\mu$.
\end{prop}

This is essentially Theorem 3.2 of \cite{KM}. The following proof is 
a transcription of the original proof in \cite{KM}, made slightly more
elementary in the sense that we avoid the use of Theorem 2.13 
(involving configurations of subspaces in spaces of dimension higher than
$3$) of
\cite{KM} at the cost of a determinant-computation.

\begin{proof} We have to prove that 
$$\hbox{lk}(L_i(t),L_j(t))=\hbox{sign}\left(\begin{array}{c}
P_\omega^i(t)-P_\alpha^i(t)\\
P_\alpha^j(t)-P_\omega^i(t)\\
P_\omega^j(t)-P_\alpha^j(t)\end{array}\right)$$
is well-defined (takes a constant value).
The matrix involved is given by
$$\left(\begin{array}{ccc}
(t-1)i-t\mu(i')&-2&(1-t)\sigma(i)-ti'\\
(1-t)j+t\mu(i')&2&(t-1)\sigma(i)+tj'\\
(t-1)j-t\mu(j')&-2&(1-t)\sigma(j)-tj'\end{array}\right)$$
and its determinant $p$ equals 
$$\begin{array}{l}
\displaystyle 2\big((i-j)(\sigma(i)-\sigma(j))+(i'-j')(\mu(i')-\mu(j'))\big)t^2
\\
\displaystyle \qquad -4(i-j)(\sigma(i)-\sigma(j))t+2(i-j)(\sigma(i)-\sigma(j))
\ .\end{array}$$
The discriminant (with respect to $t$) of $p$ given by
$$-16(i-j)(\sigma(i)-\sigma(j))(i'-j')(\mu(i')-\mu(j'))$$
is strictly negative for $i\not=j$ which shows that $p$ is non-zero
for $t\in{\mathbb R}$.
This proves that ${\mathcal L}(t)$ is a skew configuration for
$t\in \R$.

For $t=0$ we get a spindle with axis $(\R,1,0)$ and directrix $(0,-1,\R)$
realizing the spindle $\sigma$ with lines
$$L_i=s(i,1,0)+(1-s)(0,-1,\sigma(i)),\ i=1,\dots,n\ s\in\R\ .$$

For $t=1$ we get a spindle with axis $(0,1,\R)$ and directrix
$(-\R,-1,0)$ (where $-\R$ denotes the real line endowed with the opposite
order) realizing $\mu$ with lines
$$L'_{i'}=s(0,1,i')+(1-s)(-\mu(i'),-1,0),\ i'=1,\dots,n\ ,s\in \R.$$ 
\end{proof}

The proof of Theorem \ref{KM3.2} is immediate.

\begin{proof}[Proof of Corollary \ref{maincor}] Spindle-equivalent 
permutations give rise to isotopic spindles and their linking
matrices are thus switching-equivalent.

On the other hand, given two permutations $\sigma,\sigma'$
with switching-equivalent linking
matrices, Theorem \ref{KM3.2} yields an isotopy
between the associated spindle-configurations and Theorem \ref{main}
implies that $\sigma,\sigma'$ are spindle-equivalent.
\end{proof}


\section{Spindlegenus}\label{spindlegenus}

This section describes a topological invariant 
of permutations up to spindle-equivalence. 
This yields an invariant for
spindle-configurations by Theorem \ref{main}.

Let $P$ be a regular polygon with $n$ 
edges $E_1,\dots,E_n$ in clockwise cyclical order. 
Reading indices modulo $n$ we orient 
the edge $E_i$ from $E_i\cap E_{i-1}$ 
and denote by $E_{-i}$ the edge of opposite orientation.
Consider a second polygon $P'$ with edges $E'_1,\dots,E'_n$
obtained from $P$ by an orientation-reversing isometry (e.g. an orthogonal
symmetry with respect to a line). 
Given a permutation $\sigma:  \{1,\dots,n\}\longrightarrow
\{1,\dots,n\}$, gluing the oriented edge $E_i\in P$ onto the oriented 
edge $E'_{\sigma(i)}\in P'$ for $1\leq i\leq n$ yields a compact 
orientable surface $\Sigma(\sigma)$. We call the genus 
$g(\sigma)\in {\mathbb N}$ of $\Sigma(\sigma)$ the
{\it spindlegenus} of $\sigma$. 

This construction can be generalized as follows:

A {\it signed} permutation is a permutation
of the set $\{\pm 1,\dots,\pm n\}$ such that $\tilde\sigma(-i)=-\tilde
\sigma(i)$. The group of all signed permutations is the full group of 
all isometries acting on the regular $n-$dimensional 
cube $[-1,1]^n\subset \R^n$. 
Such a permutation can be graphically represented
by segments $[(i,0),(\vert\tilde\sigma(i)\vert,1)]$ carrying
signs $\epsilon_i=\frac{\vert\tilde\sigma(i)\vert}{\tilde\sigma(i)}
\in \{\pm 1\}$. 
The notion of spindle-equivalence extends
to signed permutations in the obvious way.
The construction of the compact surface $\Sigma(\sigma)$ can now be 
applied to a signed permutation $\tilde\sigma$ (glue the 
oriented edge $E_i\in P$ onto the oriented edge $E'_{\tilde\sigma(i)}\in P'$ 
for all $i=1,\dots,n$) and yields a generally 
non-orientable surface $\Sigma(\tilde\sigma)$ presented as
a $CW-$complex with two open $2-$cells (corresponding to the 
interiors of the polygons $P$ and $P'$), $n$ open $1-$cells
(the edges of $P$ or $P'$) and with
$v(\tilde \sigma)$ points or $0-$cells. The surface 
$\Sigma(\tilde \sigma)$ is orientable if and only if $\vert 
\sum_{i=1}^n \epsilon_i\vert=n$,
i.e. if $\tilde \sigma$ is either an ordinary permutation or 
the opposite of an ordinary permutation. The classification of compact 
surfaces (cf. for instance \cite{Massey}) shows
that, up to homeomorphisms, the surface $\Sigma(\tilde\sigma)$ is
completely described by $n$, $v(\tilde \sigma)$ and orientability.

\begin{prop}\label{vinvariant_prop} If $\tilde \sigma,\ \tilde
  \sigma'$ are two signed spindle-equi\-valent permutations,
then $v(\tilde\sigma)=v(\tilde \sigma')$ and the compact surfaces
$\Sigma(\tilde\sigma)$ and $\Sigma(\tilde \sigma')$ are homeomorphic.
\end{prop}

\begin{cor}\label{spindlegenus_prop}
Two permutations $\sigma$ and $\sigma'$ which are
spindle-equi\-valent have the same genus
(i.e. $g(\sigma)=g(\sigma')$).
\end{cor}

The following table lists the multiplicities for the spindlegenus
$g(\sigma)$ (related to the number $v(\sigma)$ of vertices in
the $CW-$complex considered above by the formula
$\chi(\Sigma)=2-n+v(\sigma)=2-2g(\sigma)$ for the Euler characteristic of 
$\Sigma(\sigma)$)
of permutations (normalized by $\sigma(1)=1$, multiply
by $n$ in order to get the corresponding numbers for not necessarily
normalized permutations):
$$\begin{array}{|c|c|c|c|c|c|c|}
\hline
n & g=0 & g=1 & g=2 & g=3 & g=4 & g=5 \cr
\hline
1  & 1 &     &       &         &          & \cr
2  & 1 &     &       &         &          & \cr
3  & 1 &   1 &       &         &          & \cr
4  & 1 &   5 &       &         &          & \cr
5  & 1 &  15 &     8 &         &          & \cr
6  & 1 &  35 &    84 &         &          & \cr
7  & 1 &  70 &   469 &     180 &          & \cr
8  & 1 & 126 &  1869 &    3044 &          & \cr
9  & 1 & 210 &  5985 &   26060 &     8064 & \cr
10 & 1 & 330 & 16401 &  152900 &   193248 & \cr
11 & 1 & 495 & 39963 &  696905 &  2286636 &   604800 \cr
12 & 1 & 715 & 88803 & 2641925 & 18128396 & 19056960 \cr
\hline
\end{array}$$
(see also Sequence A60593 of \cite{Sl}).

Given a permutation $\sigma$ of $\{1,\dots,n\}$ and a sequence of signs
$(\epsilon_1,\dots,\epsilon_n)\in\{\pm 1\}^n$, we consider the signed
permutation $\sigma_\epsilon:i\longmapsto \epsilon_i\sigma(i)$ and set
$$V_\sigma(t,z)=\sum_{(\epsilon_1,\dots,\epsilon_n)\in\{\pm 1\}^n}
z^{v(\sigma_\epsilon)}\ t^{\sum_i\epsilon_i}\in{\mathbf Z}[t,\frac{1}{t},z]\ .$$

\begin{cor} \label{vpolynomial}
The application $\sigma\longmapsto V_\sigma(t,z)$ is
well-defined on spindle-equivalence classes and satisfies
$$V_{\overline\sigma}(t,z)=V_\sigma(\frac{1}{t},z)$$
where $\overline\sigma(i)=n+1-\sigma(i)$ denotes a mirror spindle-permutation
of $\sigma$.
\end{cor}

\begin{proof}[Proof of Proposition \ref {vinvariant_prop}]
The number $v(\tilde \sigma)$ is obviously invariant under circular moves.

The compact surface $\Sigma(\tilde\sigma)$ is orientable
if and only if $\sum_{i=1}^n \epsilon_i\in\{\pm n\}$ and $\sum_{i=1}^n
\epsilon_i$ is preserved under spindle-equivalence.
It is thus enough to show
the invariance of the number $v(\sigma)$ of vertices in the 
$CW-$complex representing $\Sigma(\tilde
\sigma)$ under horizontal and vertical reflections.

This number $v=v(\tilde \sigma)$ can be computed graphically as follows:
Re\-present the signed permutation $\tilde\sigma(i)$ of $\sigma$ by 
drawing $n$ segments joining the
points $(i,0),\ i=1,\dots,n$ to $(\vert \sigma(i)\vert,1)$ as shown
in Figure \ref{genus_example} for the signed permutation
$\tilde \sigma(1)=-2,\ \tilde \sigma(2)=3,\ 
\tilde \sigma(3)=1$.

\begin{figure}[h]
\epsfysize=5.5cm
\epsfbox{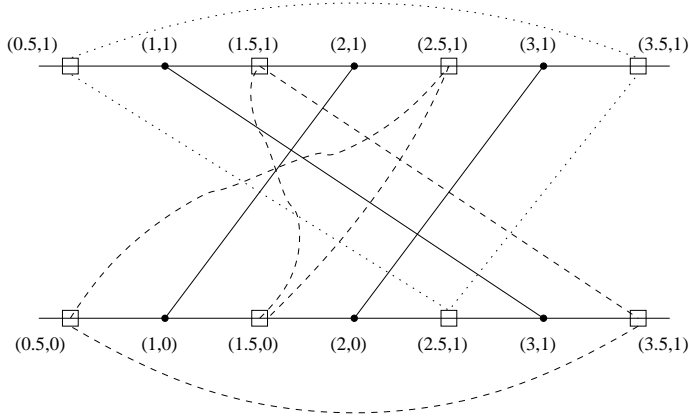}
\caption{Example for the computation of 
$v(\tilde \sigma)$}\label{genus_example}
\end{figure}                                                                 

Add $2n+2$ additional points $(0.5,0),(1.5,0),\dots,(n+0.5,0)$ and
$(0.5,1),$ $(1.5,1),\dots,(n+0.5,1)$ by drawing all points at height $y=0$ or 
$1$ which are at distance $1/2$ from an endpoint of a segment
$[(i,0),(\vert\tilde\sigma(i)\vert,1)]$. Join $(0.5,0)$,$(n+0.5,0)$
(respectively $(0.5,1)$, $(n+0.5,1)$) by dotted convex (respectively
concave) arcs and join the points $(i\pm 0.5,0)$,$(\vert
\tilde \sigma(i)\vert \pm \epsilon_i
0.5,1)$
by dotted arcs  with 
$\epsilon_i=\frac{\tilde \sigma(i)}{\vert\tilde\sigma(i)\vert}$ 
corresponding to the sign of the $i-$th segment. This yields 
a graph with vertices $(0.5,0),(1.5,0),\dots,(n+0.5,0),
(0.5,1),$ $(1.5,1),\dots,(n+0.5,1)$ of degree $2$ 
by considering all dotted arcs as edges. The number of
connected components of this graph ($2$ in Figure \ref{genus_example})
equals $v(\tilde\sigma)$. This can be seen as follows: The interior
of the polygons $P,P'$ correspond to the half-planes
$y<0$ and $y>1$. Vertices of $P,P'$ correspond to the points
$(0.5,0),(1.5,0),\dots,(n+0.5,0),
(0.5,1),$ $(1.5,1),\dots,(n+0.5,1)$ where the pairs of points
$(0.5,0),(n.5,0)$ and $(0.5,1),(n.5,1)$ have to be glued together
(achieved by additional arcs joining these points).
The segments $(i,0),(\vert\tilde \sigma(i)\vert,1)$ represent glued
edges with dotted arcs joining vertices identified under gluing.

Let us consider the local situation around a block $B$ of $\sigma$.
The (internal) dotted arcs associated to lines of $B$ 
connect the four boundary points adjacent to $B$
in one of the three ways depicted in Figure
\ref{genus_proof}. The proof is now obvious since each of 
these three situations is invariant under vertical and horizontal reflections. 
\end{proof}

\begin{figure}[h]
\epsfysize=5cm
\epsfbox{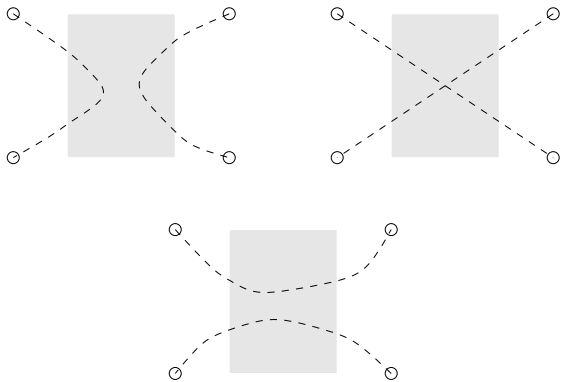}
\caption{Three local situations around a block}\label{genus_proof}
\end{figure}                                                                 

We leave the easy proofs of Corollaries \ref{spindlegenus_prop}
and \ref{vpolynomial} to the reader.

\begin{rem} \label{spindle_generalizations}
Many similar invariants of spindle-permutations, 
up to spindle-equi\-valence,
can be defined similarly by considering a set ${\mathcal S}$
of subsets of lines in $\sigma$
which is defined in a topological way (e.g. triplets of lines of given isotopy
class, subsets of $k$ lines defining a subspindle in a set of prescribed
spindle-equivalence classes, or a given subset of lines arising in the 
Euler partition of a spindle having  $2n$ lines) and by considering
all sign sequences $(\epsilon_1,\dots,\epsilon_n)\in\{\pm 1\}^n$
such that $\{i\ \vert\ \epsilon_i=-1\}\in {\mathcal S}$. The corresponding
sum $\sum z^{v(\sigma_\epsilon)}$ is then of course well-defined for 
spindles and thus for spindle-permutations, up to equivalence.
\end{rem}

We end this section with the following natural questions: 

\begin{enumerate}
\item By Corollary \ref{maincor}, the 
(complete) spindlegenus factorize through switching classes
associated to spindles. Can
the spindlegenus (or some related invariant)
be extended to all switching classes?

\item Can the (complete) spindlegenus
be extended to skew configurations which are not 
spindles?
\end{enumerate}

\section{Spindle structures for switching classes}\label{algo}

The existence
of a linking matrix of a spindle in a given switching class is
a natural question which we want to address algorithmically 
in this section.

The following algorithm exhibits a spindle-permutation with 
linking matrix in a given switching class or proves non-existence
of such a permutation. By Corollary \ref{maincor}, such a spindle-permutation
is unique, up to spindle-equivalence.

\medskip
\begin{algo} \ 
\begin{itemize}
\item[{\bf Initial data.}] A natural number $n$ and a switching class
represented by a symmetric matrix
$X$ of order $n$ with rows and columns indexed by $\{1,\dots,n\}$
and coefficients $x_{i,j}$ satisfying

$$\begin{array}{ll}
x_{i,i}=0,\qquad &1\leq i\leq n\ ,\cr
x_{i,j}=x_{j,i}\in \{\pm 1\}\ ,\quad &1\leq i\not= j\leq n\ .
\end{array}$$

\item[{\bf Initialization.}] Conjugate the symmetric matrix $X$ by the
diagonal matrix with diagonal coefficients $(1,x_{1,2},x_{1,3},
\dots,x_{1,n})$. Set $\gamma(1)=\gamma(2)=1,\ \sigma(1)=1$ and $k=2$.

\item[{\bf Main loop.}] 

Replace $\gamma(k)$ by $\gamma(k)+1$ and set
$$\sigma(k)=1+\sharp\{j\ \vert\ x_{\gamma(k),j}=-1\}+\sum_{s=1}^{k-1}
x_{\gamma(s),\gamma(k)}\ .$$

Check the following conditions:
\begin{enumerate}
\item $\gamma(k)\not=\gamma(s)$ for $s=1,\dots,k-1$.
\item $x_{\gamma(k),\gamma(s)}=\hbox{sign}(\sigma(k)-\sigma(s))$ for
$s=1,\dots,k-1$
(where $\hbox{sign}(0)=0$ and $\hbox{sign}(x)=\frac{x}{\vert x\vert}$
for $x\not=0$).
\item for $j\in \{1,\dots,n\}\setminus\{\gamma(1),\dots,\gamma(k)\}$ and
$1\leq s<k$:\\
if $x_{j,\gamma(s)}\ x_{\gamma(s),\gamma(k)}=-1$ then 
$x_{j,\gamma(k)}=x_{j,\gamma(s)}$.
\end{enumerate}

\medskip

\noindent
If all conditions are fulfilled then:
\begin{itemize}
\item[] if $k=n$ print all the data (mainly the spindle-permutation 
$i\longmapsto \sigma(i)$ and perhaps also the conjugating permutation
$i\longmapsto \gamma(i)$) and stop.
\item[] if $k<n$ then set $\gamma(k+1)=1$, replace $k$ by $k+1$ and redo the
main loop.
\end{itemize}

\medskip

\noindent
If at least one of the above conditions is not fulfilled then:
\begin{itemize}

\item[] while $\gamma(k)=n$ replace $k$ by $k-1$.

\item[] if $k=1$: print \lq\lq no spindle structure exists for
this switching class'' and stop.

\item[] if $k>1$: redo the main loop.

\end{itemize}

\end{itemize}
\end{algo}

\medskip

\noindent
{\bf Explanation of the algorithm.}  The {\bf initialization} is
actually a
normalization: we assume that the first row of the matrix represents the
first line of a spindle-permutation $\sigma$ normalized to $\sigma(1)=1$
(up to a circular move, this can always be assumed for a 
spindle-permutation).

The {\bf main loop} assumes that row number $\gamma(k)$ of $X$ contains
the crossing data of the $k-$th line $L_k$
(assuming that the rows representing the crossing data of $L_1,\dots,
L_{k-1}$ are correctly chosen). 
The image $\sigma(k)$ of $k$ under a spindle-permutation is 
then uniquely defined and given by the formula used in the main loop. 

One has to check three necessary conditions:
\begin{itemize}
\item The first condition checks that row number $\gamma(k)$ has
not been used before. 
\item The second condition checks the consistency of the choice for $\gamma(k)$
with all previous choices.
\item If the third condition is violated, then the given choice of rows
$\gamma(1),\dots,\gamma(k)$ cannot be extended up to $k=n$.
\end{itemize}

The algorithm runs correctly even without checking for Condition (3). 
However, it looses much of its interest: Condition (3) is very
strong (especially in the case of non-existence of a spindle
structure) and ensures a fast running time.

The algorithm, in the case of success, produces two permutations
$\sigma$ and $\gamma$. The linking matrix of the spindle
permutation $\sigma$ is in the switching class as $X$ and 
$\gamma$ yields a conjugation between these two matrices.
More precisely:
$$x_{\gamma(i),\gamma(j)}=\hbox{sign}((i-j)(\sigma(i)-\sigma(j)))$$
under the assumption $x_{1,i}=x_{i,1}=1$ for $2\leq i\leq n$.

Failure of the algorithm (the algorithm stops
after printing ``no spindle structure exists for this switching 
class'') proves non-existence of a spindle structure in
the switching class of $X$.

In practice, the average running time of this algorithm should be 
of order $O(n^3)$ or perhaps $O(n^4)$. 
Indeed Condition (3) is only
very rarely satisfied for a wrong choice of $\gamma(k)$ with $k>2$.
On the other hand, for a switching class containing a crossing
matrix of a spindle, checking all cases of Condition (3) needs at
least $O(n^3)$ operations (or more precisely ${n-1\choose 3}$
operations after suppressing the useless comparisons involving
$\gamma(1)=1$).

\begin{rem} The algorithm can be improved. Condition (3) can be
made considerably stronger.
\end{rem}

\section{Computational results and lower bounds for the number of
  non-isotopic configurations}\label{enumeration}

In this section, we describe some computational results.

The number of skew configurations,up to isotopy, having $n\leq 7$
lines are known (see the survey of Viro and
Drobotukhina \cite{VD} and the results of Borobia and Mazurovskii
\cite{BM1}, \cite{BM2}, \cite{Ma}): 

\medskip

$$
\begin{array}{|c|c|}
\hline
{\rm Lines} & {\rm Isotopy\ classes}  \\
\hline
2 & 1  \\
3 & 2  \\
4 & 3  \\
5 & 7  \\
6 & 19 \\
7 & 74 \\
\hline
\end{array}
$$

\medskip

The following table enumerates the number of switching classes of
order $6-9$. The middle row shows the number of distinct polynomials
which arise as characteristic polynomials of switching classes (this
is of course the same as the number of conjugacy classes under the
orthogonal group $O(n)$ of matrices representing switching classes).

\medskip

$$
\begin{array}{|c|c|c|}
\hline
{\rm Lines} & {\rm Characteristic} & {\rm Switching} \\
            & {\rm polynomials}    & {\rm classes} \\
\hline
6 & 16 & 16 \\
7 & 54 & 54 \\
8 & 235 & 243 \\
9 & 1824 & 2038\\
\hline
\end{array}
$$

\medskip

In fact, one can use 
representation theory of the symmetric groups in order to
derive a formula for the number of switching classes of given order
(see \cite{MS} and Sequence A2854 in \cite{Sl}).

\medskip

The map
$$\{\hbox{configurations of skew
  lines}\}\longrightarrow \{ \hbox{switching classes} \}$$
is perhaps not surjective in general (there seems to be an unpublished  
counterexample of Peter Shor for $n=71$, see \cite[Section 3]{CP}). 
We rechecked however a claim of Crapo and Penne (Theorem 5 of Section 4)
stating that all $243$ switching classes of order $8$ 
arise as linking matrices of suitable skew configurations of $8$
lines (the corresponding result holds also for fewer lines). 
There are thus at least 243 isotopy classes of configurations containing 
$8$ skew lines, 180 of them are spindle classes.

\medskip

The following table shows the number of spindle-permutations, 
up to equivalence, for $n\leq 13$. We also indicate the number of
amphicheiral spindle-permutations, up to equivalence.
  
\medskip

$$
\begin{array}{|c|c|c|}
\hline
n & \hbox{spindle classes}& \hbox{ amphicheiral classes}   \\
\hline
1  & 1       & 1  \\
2  & 1       & 1   \\
3  & 2       & 0   \\
4  & 3       & 1   \\
5  & 7       & 1  \\
6  & 15      & 3   \\
7  & 48      & 0   \\
8  & 180     & 12    \\
9  & 985     & 5   \\
10 & 6867    & 83   \\
11 & 60108   & 0     \\
12 & 609112  & 808   \\
13 & 6909017 & 47   \\
\hline
\end{array}
$$

\medskip

Assertion (ii) of Proposition \ref{mirror} 
explains of course the non-existence of amphicheiral classes 
for $n\equiv 3\pmod 4$.

\section*{Acknowledgments} 
We express our thanks to two anonymous referees of
a previous version (consisting mainly of a flawed proof of Theorem \ref{main})
for helpful comments and remarks. The first author thanks also 
D. Matei for his interest and valuable comments.

The second author wish to thank Institut Fourier and Mikhail
Zaidenberg for hosting his stay.

\end{document}